\input amssym.def
\input amssym.tex

\def\item#1{\vskip1.3pt\hang\textindent {\rm #1}}


\tolerance=300
\pretolerance=200
\hfuzz=1pt
\vfuzz=1pt


\magnification=\magstep0


\hoffset=0.6in
\voffset=0.8in

\hsize=5.8 true in 
\vsize=8.5 true in
\parindent=25pt
\mathsurround=1pt
\parskip=1pt plus .25pt minus .25pt
\normallineskiplimit=.99pt

\countdef\revised=100
\mathchardef\emptyset="001F 
\chardef\ss="19
\def\3{\ss}
\def\anf{$\lower1.2ex\hbox{"}$}
\def\frac#1#2{{#1 \over #2}}
\def\>{>\!\!>}
\def\<{<\!\!<}

\def\ssarr{\hbox to 30pt{\rightarrowfill}}
\def\sarr{\hbox to 40pt{\rightarrowfill}}
\def\arr{\hbox to 60pt{\rightarrowfill}}

\def\larr{\hbox to 60pt{\leftarrowfill}}
\def\Arr{\hbox to 80pt{\rightarrowfill}}

{}

\def\Str{\mathop{\rm Str}\nolimits}

\def\ad{\mathop{\rm ad}\nolimits}

\def\Asym{\mathop{\rm Asym}\nolimits}

\def\Aut{\mathop{\rm Aut}\nolimits}

\def\det{\mathop{\rm det}\nolimits}

\def\End{\mathop{\rm End}\nolimits}

\def\Gl{\mathop{\rm Gl}\nolimits}

\def\Hom{\mathop{\rm Hom}\nolimits}%
\def\id{\mathop{\rm id}\nolimits} 
\def\im{\mathop{\rm im}\nolimits}

\def\OO{{\rm O}}



\def\Sym{\mathop{\rm Sym}\nolimits}

\def\Sp{\mathop{\rm Sp}\nolimits}


\def\UU{{\rm U}}

\def\OO{{\rm O}}

\def\pr{\mathop{\rm pr}\nolimits}

\def\0{{\bf 0}}
\def\1{{\bf 1}}

\def\a{{\frak a}}

\def\b{{\frak b}}

\def\e{{\frak e}}
\def\f{{\frak f}}
\def\g{{\frak g}}

\def\h{{\frak h}}

\def\k{{\frak k}}

\def\m{{\frak m}}

\def\n{{\frak n}}

\def\p{{\frak p}}
\def\q{{\frak q}}

\def\s{{\frak s}}

\def\sL{{\frak {sl}}}

\def\X{{\frak X}}

\def\K{{\Bbb K}} 
\def\LL{{\Bbb L}} 
 
\def\N{{\Bbb N}} 
 
\def\P{{\Bbb P}}

\def\Z{{\Bbb Z}} 

\def\:{\colon}  
\def\.{{\cdot}}
\def\|{\Vert}
\def\bsk{\bigskip}

\def\giantskip{\vskip2\bigskipamount}
\def\gsk{\giantskip}

\def\msk{\medskip}

\def\ssk{\smallskip}

\def\bbr{\bigbreak}
\def\giantbreak{\par \ifdim\lastskip<2\bigskipamount \removelastskip
         \penalty-400 \giantskip\fi}

\def\nin{\noindent}
\def\cen{\centerline}
\def\pagebreak{\vskip 0pt plus 0.0001fil\break}
\def\linebreak{\break}

\def\epsilon{\varepsilon}

\def\nin{\noindent}

\def\pder#1,#2,#3 { {\partial #1 \over \partial #2}(#3)}
\def\pde#1,#2 { {\partial #1 \over \partial #2}}
\def\phi{\varphi}


\def\tilde{\widetilde}

\font\eightrm=cmr8


\font\bfone=cmbx10 scaled\magstep1 
\font\bftwo=cmbx10 scaled\magstep2 

\def\qed{{\unskip\nobreak\hfil\penalty50\hskip .001pt \hbox{}\nobreak\hfil
          \vrule height 1.2ex width 1.1ex depth -.1ex
           \parfillskip=0pt\finalhyphendemerits=0\medbreak}\rm}

\def\qeddis{\eqno{\vrule height 1.2ex width 1.1ex depth -.1ex} $$
                   \medbreak\rm}

\def\Lemma #1. {\bigbreak\vskip-\parskip\noindent{\bf Lemma #1.}\quad\it}

\def\Sublemma #1. {\bigbreak\vskip-\parskip\noindent{\bf Sublemma #1.}\quad\it}

\def\Proposition #1. {\bigbreak\vskip-\parskip\noindent{\bf Proposition #1.}
\quad\it}

\def\Corollary #1. {\bigbreak\vskip-\parskip\nin{\bf Corollary #1.}
\quad\it}

\def\Theorem #1. {\bigbreak\vskip-\parskip\noindent{\bf Theorem #1.}
\quad\it}

\def\Definition #1. {\rm\bigbreak\vskip-\parskip\noindent{\bf Definition #1.}
\quad}

\def\Remark #1. {\rm\bigbreak\vskip-\parskip\noindent{\bf Remark #1.}\quad}

\def\Example #1. {\rm\bigbreak\vskip-\parskip\noindent{\bf Example #1.}\quad}

\def\Problems #1. {\bigbreak\vskip-\parskip\noindent{\bf Problems #1.}\quad}
\def\Problem #1. {\bigbreak\vskip-\parskip\noindent{\bf Problem #1.}\quad}

\def\Conjecture #1. {\bigbreak\vskip-\parskip\noindent{\bf Conjecture #1.}\quad}

\def\Proof#1.{\rm\par\ifdim\lastskip<\bigskipamount\removelastskip\fi\smallskip
            \noindent {\bf Proof.}\quad}

\def\Axiom #1. {\bigbreak\vskip-\parskip\noindent{\bf Axiom #1.}\quad\it}

\def\Satz #1. {\bigbreak\vskip-\parskip\noindent{\bf Satz #1.}\quad\it}

\def\Korollar #1. {\bbr\vskip-\parskip\nin{\bf Korollar #1.} \quad\it}

\def\Bemerkung #1. {\rm\bigbreak\vskip-\parskip\noindent{\bf Bemerkung #1.}
\quad}

\def\Beispiel #1. {\rm\bigbreak\vskip-\parskip\noindent{\bf Beispiel #1.}\quad}
\def\Aufgabe #1. {\rm\bigbreak\vskip-\parskip\noindent{\bf Aufgabe #1.}\quad}

\def\Beweis#1. {\rm\par\ifdim\lastskip<\bigskipamount\removelastskip\fi
           \smallskip\noindent {\bf Beweis.}\quad}

\nopagenumbers

\def\date{ }

\def\title{Title ??}
\def\author{Author ??}

\def\thanks#1{\footnote*{\eightrm#1}}

\def\rightheadline{\hfil{\eightrm\author}\hfil\tenbf\folio}
\def\leftheadline{\tenbf\folio
\hfil}
\headline={\vbox{\line{\ifodd\pageno\rightheadline\else\leftheadline\fi}}}

\def\firstheadline{}
\def\firstfootline{\cen{\rm\folio}}

\def\seite #1 {\pageno #1
               \headline={\ifnum\pageno=#1 \firstheadline
               \else\ifodd\pageno\rightheadline\else\leftheadline\fi\fi}
               \footline={\ifnum\pageno=#1 \firstfootline\else{}\fi}}

\newdimen\dimenone
 \def\checkleftspace#1#2#3#4{
 \dimenone=\pagetotal
 \advance\dimenone by -\pageshrink   
 \ifdim\dimenone>\pagegoal          
   \else\dimenone=\pagetotal
        \advance\dimenone by \pagestretch
        \ifdim\dimenone<\pagegoal
          \dimenone=\pagetotal
          \advance\dimenone by#1         
          \setbox0=\vbox{#2\parskip=0pt                
                     \hyphenpenalty=10000
                     \rightskip=0pt plus 5em
                     \noindent#3 \vskip#4}    
        \advance\dimenone by\ht0
        \advance\dimenone by 3\baselineskip   
        \ifdim\dimenone>\pagegoal\vfill\eject\fi
          \else\eject\fi\fi}


\def\sectionheadline #1{\bigbreak\vskip-\lastskip
      \checkleftspace{1.1cm}{\bf}{#1}{\bigskipamount}
         \vbox{\vskip1.1cm}\cen{\bfone #1}\bsk}

\def\lsectionheadline #1 #2{\bigbreak\vskip-\lastskip
      \checkleftspace{1.1cm}{\bf}{#1}{\bigskipamount}
         \vbox{\vskip1.1cm}\cen{\bfone #1}\msk \cen{\bfone #2}\bsk}

\def\lchapterheadline #1 #2{\bigbreak\vskip-\lastskip\indent\vskip3cm
                       \cen{\bftwo #1} \msk \cen{\bftwo #2} \gsk}
\def\llsectionheadline #1 #2 #3{\bigbreak\vskip-\lastskip\indent\vskip1.8cm
\cen{\bfone #1} \msk \cen{\bfone #2} \msk \cen{\bfone #3} \nobreak\bsk\nobreak}


\newtoks\literat
\def\[#1 #2\par{\literat={#2\unskip.}%
\hbox{\vtop{\hsize=.15\hsize\nin [#1]\hfill}
\vtop{\hsize=.82\hsize\nin\the\literat}}\par
\vskip.3\baselineskip}

\def\references{
\sectionheadline{\bf References}
\frenchspacing

\entries\par}

\mathchardef\emptyset="001F 
\def\address{Author: \tt$\backslash$def$\backslash$address$\{$??$\}$}

\def\firstpage{\nin
{\obeylines \parindent 0pt }
\vskip2cm
\centerline{\bfone\title}
\gsk
\centerline{\bf\author}
\vskip1.5cm \rm}

\def\addresstwo{}

\def\dlastpage{\par\vbox{\vskip1cm\nin
\line{
\vtop{\hsize=.5\hsize{\parindent=0pt\baselineskip=10pt\nin\address}}
\quad 
\vtop{\hsize=.42\hsize\nin{\parindent=0pt
\baselineskip=10pt\addresstwo}}
\hfill} }}


\pageno=1
\def\title
{Inner ideals and intrinsic subspaces in linear pair geometries}

\def\author{Wolfgang Bertram, Harald L\"owe}

\def\date{13.06.2006}
\def\leftheadline{\tenbf\folio\hfil{\eightrm\title}\hfil\eightrm\date}
\def\address{Harald L\"owe

Technische Universit\"at Braunschweig

Institut Computational Mathematics

Pockelsstrasse 14

D - 38106 Braunschweig

e-mail: h.loewe@tu-bs.de}

\def\addresstwo{Wolfgang Bertram

Universit\'e Henri Poincar\'e (Nancy I) -- Institut Elie Cartan

B.P. 239

F - 54506 Vand\oe uvre-l\`es-Nancy Cedex 

e-mail: bertram@iecn.u-nancy.fr}


\def\u{{\frak u}}
\def\v{{\frak v}}
\def\b{{\frak b}}

\def\o{{\frak o}}
\def\k{{\frak k}}


\def\YY{{\cal Y}}
\def\XX{{\cal X}}
\def\FF{{\cal F}}
\def\LL{{\cal L}}

\def\WW{{\cal W}}
\def\II{{\cal I}}
\def\JJ{{\cal J}}

\def\Lag{{\cal L}ag}
\def\Gras{{\cal G}ras}

\def\PE{\mathop{\rm PE}\nolimits}

\def\Gl{\mathop{\rm GL}\nolimits}

\def\Str{\mathop{\rm Str}\nolimits}
\def\pr{\mathop{\rm pr}\nolimits}

\def\Asym{\mathop{\rm Asym}\nolimits}

\def\S{{\bf \Sigma}}
\def\P{{\bf \Pi}}

\firstpage

\bigskip
\centerline{Wolfgang Bertram }
\centerline{Universit\'e Henri Poincar\'e (Nancy I) -- Institut Elie Cartan}
\centerline{B.P. 239}
\centerline{F - 54506 Vand\oe uvre-l\`es-Nancy Cedex} 
\centerline{e-mail: bertram@iecn.u-nancy.fr }

\msk
\centerline{Harald L\"owe}
\centerline{Technische Universit\"at Braunschweig --
Institut Computational Mathematics}
\centerline{Pockelsstrasse 14}
\centerline{D - 38106 Braunschweig}
\centerline{e-mail: h.loewe@tu-bs.de}

\bigskip
\nin {\bf Abstract.} 
We introduce the notion of {\it intrinsic subspaces} of
{\it linear} and {\it affine pair geometries}, which generalizes
the one of projective subspaces of projective spaces. 
We prove that, when the affine pair geometry is the {\it projective
geometry of a Lie algebra} introduced in [BeNe04], such intrinsic
subspaces correspond to {\it inner ideals} in the associated Jordan
pair, and we investigate the case of intrinsic subspaces defined
by the Peirce-decomposition which is related to $5$-gradings of the
projective Lie algebra. These examples, as well as the examples of
general and Lagrangian flag geometries, lead to the 
conjecture that geometries of
intrinsic subspaces tend to be themselves linear pair geometries.

\bigskip
\nin {\bf Contents:}

\ssk
\nin I. The axiomatic approach

1. Linear pair geometries

2. Intrinsic subspaces in linear pair geometries

\nin II. The main examples

3. Flag geometries

4. Lagrangian flag geometries

5. Geometries of filtrations in Lie algebras

\nin Appendix A: Inner ideals in matrix spaces

\nin Appendix B: Inner ideals and Jordan theory

\bigskip

\nin {\bf AMS-classification:}
17C37,  
17C27, 
51N30, 
51B05 

\bigskip

\nin {\bf Key words:} linear and affine pair geometries, intrinsic
subspace, inner ideal, Jordan pair, graded and filtered Lie algebras,
flag geometries

\vfill \eject

\sectionheadline{Introduction}

The notion of {\it subspace of a projective space}
is classical, and these subspaces have the nice property that,
when looking at them from the affine point of view (i.e, after having 
removed some hyperplane at infinity) they just look like usual affine
subspaces. Since projective spaces are generalized by
{\it Grassmannians} (say, the Grassmannian $\Gras_{p,q}(\K)$
of $p$-spaces in some $p+q$ dimensional
$\K$-vector space $W$), one may ask for a natural notion of subspace
also in this and similar cases.
One may expect, then, that the collection of all such subspaces forms
a new kind of ``geometry'', with probably interesting properties and
interacting with the theory of the Grassmannian geometry we started with.
In this work, we will define such a kind of subspaces,
 called {\it intrinsic subspaces}, in a rather general framework 
which includes Grassmannians and many other examples (finite or
infinite dimensional, or more generally in geometries over rings); 
we will give
explicit constructions of such subspaces for many cases, and relate them
to well-known objects that are called, in the theory of Jordan algebraic
structures, {\it inner ideals} 
 and which have  attracted interest
 as possible material to build on them, e.\ g., quantum mechanical
theories with a geometrical flavor (see, e.g., work of J.R.\ Faulkner [F80] 
and Chapter 9 of [BeNe05] for a link with the present work). 
Moreover, this approach permits to give a geometric definition of
the notion of
{\it rank} (or, using a terminology introduced by L.\ K.\ Hua [Hua45]
for the classical matrix spaces, of {\it arithmetic distance}) in such 
spaces; the stratification of tangent spaces by the rank is very closely
related to the {\it generalized conformal structures} introduced by
S.\ Gindikin and S.\ Kaneyuki ([GiKa98]; cf.\ [Be00, Ch.\ IX]).
In a sense, the present approach can be seen as a rather far-reaching
generalization of these notions.

\ssk
The basic idea is easily explained at the example of the Grassmannian
geometry $\XX = \Gras_{p,q}(\K)$.
Just like the projective space, the Grassmannian is covered by ``affine
parts'', namely by the sets $\v_\alpha$ of $p$-spaces that are
complementary to some given $q$-space $\alpha$.
Thus the affine parts are naturally parametrized by the ``dual
Grassmannian'' $\XX' = \Gras_{q,p}(\K)$.
Now, in the projective case (i.e., $p=1$) it is true that a subset
which looks affinely with respect to {\it one} affinization looks
affinely with respect to {\it all} affinizations. 
However, in the general case
(i.e.\ $p>1$ and $q>1$) this is false!
On the other hand, there exist some subsets which do have this property,
and this defines precisely our intrinsic subspaces:
an {\it intrinsic subspace of $\XX=\Gras_{p,q}(\K)$} is a subset 
$\II \subset \XX$ such that all intersections
$\II_\alpha = \II \cap \v_\alpha$, where
$\alpha$ runs through $\XX'$, are affine subspaces of $\v_\alpha$.

\ssk
The question arises then: how do such intrinsic subspaces look like,
say, with respect to one fixed affinization $\v=\v_{o'}$?
Of course, they are linear (if they contain some fixed origin $o$
of $\v$); but they must satisfy some additional condition.
This additional condition is precisely the one of being an
{\it inner ideal} of $\v$.
Let us recall, for readers not familiar with Jordan theory, that
the notion of {\it inner ideal} is a simultaneous generalization
of those of left and right ideals in an associative algebra.
For instance, in the matrix algebra $A=M(n,n;\K)$ over a field $\K$,
all left ideals are, with respect to a suitable basis, described as
 sets of matrices
of the form $\bigl( {* \atop *}{0 \atop 0} \bigr)$; identifying $A$ with
$\End_\K(U)$ where $U=\K^n$, we can write them as
$I_E = \{ f \in A | \, E \subset \ker f \}$ for some fixed subspace
$E \subset U$. Simarly, right ideals are of the form
$\bigl( {* \atop 0}{* \atop 0} \bigr)$, resp.\
$J_F = \{ f \in A | \, \im f \subset F \}$ for a subspace $F \subset U$.
An {\it inner ideal} in $A = M(n,n;\K)$ is a subspace $I$ which is stable
under multiplication from the inside in the trilinear product
$$
T:A \times A \times A \to A, \quad (f,g,h) \mapsto fgh + hgf,
$$
i.e., it satisfies the condition
$T(I,A,I) \subset I$.
It is immediate that all left or right ideals of $A$ are inner, and that
an arbitrary intersection of inner ideals is again an inner ideal.
Hence  intersections of left- and right associative ideals
are again inner ideals. These intersections are matrix
spaces of the form $\bigl( {* \atop 0}{0 \atop 0} \bigr)$; more precisely,
they are described as
$$
I_{E,F}:= \{ f \in \End(U) | \, \ker(f) \subset E, \, \im(f) \subset F \}
$$
The notion of inner ideal depends only on the symmetrized product
$x \bullet y = {xy + yx \over 2}$ in the associative algebra $A$ and makes
in fact sense for any Jordan algebraic structure (-algebra, -triple system
or -pair). For instance, the algebra of matrices may be replaced by
two spaces of rectangular matrices,
$$
V^+ = M(p,q;\K) = \Hom_\K(U,W), \quad
V^- = M(q,p;\K) = \Hom_\K(W,U),
$$
where the trilinear products $T^\pm:V^\pm \times V^\mp \times
V^\pm \to V^\pm \times$ and inner ideals (in $V^+$) are defined by
the same formulae as above.
Now, coming back to the Grassmannian example, the affine part
$\v$ is naturally identified with the space $V^+$, and we show then
(Theorem 3.11) that the affine pictures $I=\v \cap \II$ of intrinsic
subspaces $\II$ are precisely the inner ideals of $V^+$.
In this case, the inner ideals can be classified:
they are all of the form $I_{E,F}$ as defined above (Proposition A.4).
Therefore the spaces $I_{E,F}$ are precisely the affine pictures of
intrinsic subspaces.
For the case of the matrix spaces, we have gathered the basic facts
on inner ideals in the first appendix (Appendix A), which uses only
elementary linear algebra. The general Jordan theoretic notions are
given in the second appendix (Appendix B), which is needed only
for Chapter 5. We hope that in this way the present work
becomes accessible to readers who are interested in geometry
but are not necessarily  specialists in Jordan theory.

\ssk
The general framework in which we can define  intrinsic subspaces,
called {\it linear pair geometries}, is defined in Chapter 1: 
it is given by a pair $(\XX^+,\XX^-)$ of sets (we think of $\XX^-$ as
the ``dual" of $\XX^+$ and therefore write also $(\XX,\XX')$), together
with a binary relation $\top$ on $\XX^+ \times \XX^-$, called
{\it transversality} or {\it remoteness} (in the Grassmann case,
$x \top \alpha$ means that $\alpha$ is a complementary subspace of $x$)
such that for every $x \in \XX^+$ there is $\alpha \in \XX^-$ with
$x \top \alpha$, and vice versa; then we assume that, for any
transversal pair $(y,\alpha)$, the set $\v_\alpha := \{ x \in \XX^+
\vert \, x \top \alpha \}$, with origin $y$, is equipped with a linear
(i.e.\ $\K$-module) structure, and vice versa.
As mentioned above, this is true in the Grassmannian case; but in that
case, moreover, the underlying affine structure of $\v_\alpha$ does
not depend on $x$; we then say that the geometry is an {\it affine
pair geometry}. Affine pair geometries have been introduced in [Be02],
where it has been shown that a special class of such geometries, called
{\it generalized projective geometries}, are essentially equivalent to
Jordan algebraic structures. We prove in this work (Theorem 5.6) that,
under this equivalence, intrinsic subspaces of the geometry correspond
to inner ideals of the Jordan structure.

\ssk
Let us now comment on the question whether the collection of all intrinsic
subspaces forms some interesting kind of geometry.
Although at present we do not have a general theory allowing to
understand these geometries, the
examples indicate that this is indeed the
case: namely, it seems that geometries of intrinsic subspaces in
linear pair geometries tend to be themselves again linear pair
geometries, but in general these geometries are no longer affine
pair geometries, even if the geometry we started with was one.
In fact, this is our main motivation to work, from the beginning on, in
this bigger class of geometries.
For instance, if we start with a Grassmann
geometry $(\XX,\XX')=(\Gras_{p,q}(\K),\Gras_{q,p}(\K))$ 
(which is an affine pair geometry), then the result stated above 
can be rephrased by saying that intrinsic subspaces
are parametrized by ``short flags'' 
$0\subset \h \subset \k \subset \K^{p+q}$, subject to the
condition $\dim \h \leq p \leq \dim \k$:
to such a flag corresponds the intrinsic subspace
$$
\II_{\h,\k} = \{ x \in \XX | \, \h \subset x \subset \k \}
$$
(such subspaces of Grassmannians are also called {\it Grassmann lines},
cf.\ [Bue95, p.\ 49]).
Thus the geometry of intrinsic subspaces corresponds, in this case,
to flag varieties of short flags. We show in Chapter 3 that
flag varieties of general linear groups (and, similarly,
flag varieties of orthogonal and symplectic groups, see Chapter 4 ) 
are indeed examples of linear, non affine pair geometries
(Theorem 3.5) 
and we construct standard examples of intrinsic subspaces for such
geometries (Theorem 3.8).
Similar observations can be made for the case of the projective geometry 
of a Lie algebra, i.e., for the  affine pair geometry of inner
$3$-filtrations of a Lie algebra $\g$ (defined in [BeNe04]):
in this case, intrinsic subspaces correspond exactly to inner ideals
in the associated Jordan pair (Theorem 5.6), and the important class
of intrinsic subspaces corresponding to inner ideals defined by
Peirce-decompositions is closely related to $5$-gradings and $5$-filtrations
of $\g$ (Theorem 5.8). 
The geometry of $5$-gradings of a Lie algebra is not yet fully understood;
we conjecture that it is also a linear pair geometry. 
Thus, at this point, several problems and topics for further research
arise naturally:

\ssk
\item{(1)} What kinds of graded Lie algebras give rise to linear
pair geometries? For $3$-graded Lie algebras, a result of this type
has been proved in [BeNe04] (see the summary in Section 5); is an
analogue true for other kinds of gradings?
\item{(2)} Do linear pair geometries of the kind mentioned in (1)
correspond to some kind of ``non-commutative Jordan structures" in
a similar way as generalized projective geometries correspond to usual
(commutative) Jordan structures?
\item{(3)} If the answer to (1) or (2) is positive, then it seems possible
to develop the theory of these geometries, similarly as for the Jordan
case in [BeNe04] and [BeNe05], in arbitrary dimension and over general base 
fields and rings, in particular furnishing interesting classes
of infinite dimensional manifolds.
\item{(4)} Finally, if it is true that geometries of intrinsic subspaces
tend to be linear pair geometries themselves, one may ask for a conceptual
proof of this fact. As a first step, it seems to be important to find a good
geometric analog of the complementation of inner ideals in Jordan pairs
(see [LoNe94]) in order to define a structure of pair geometry on the
geometry of intrinsic subspaces.

\bigskip \nin
{\bf Notation and terminology.} Throughout the paper, let $\K$ denote
a ring with unit $1$ such that $2$ is invertible in $\K$.
In some parts of the text, we will assume that $\K$ is a field.
We also use the terms {\it linear (sub)space} instead of
{\it $\K$-(sub)module} and try to use systematically gothic letters
$\v,\o,\e,\f,\ldots$ for linear spaces (or flags of linear spaces)
and calligraphic letters $\XX,\II,\FF,\LL,\ldots$ for non-linear objects. 

\bigskip \nin
{\bf Acknowledgements.}
The second author thanks the Institut Elie Cartan for hospitality and
support during his research visit in 2003, and the first author thanks
the Mathematical Institute of the TU Braunschweig for hospality 
during his visit in 2004.

\sectionheadline
{1. Linear pair geometries}

\nin
{\bf 1.1. Pair geometries.}
A {\it pair geometry} is given by a pair of sets $(\XX,\XX')$, together with 
a subset $(\XX \times \XX')^\top \subset \XX \times \XX'$, such that the following conditions are satisfied:

\ssk
\item{--} for all $\alpha \in \XX'$ there exists $x \in \XX$ such
that $(x,\alpha)\in (\XX \times \XX')^\top$, and
\item{--} for all $x \in \XX$, there exists
$\alpha \in \XX'$  such that $(x,\alpha)\in (\XX \times \XX')^\top$.

\ssk
\nin 
If $(x,\alpha)$ belongs to $(\XX \times \XX')^\top$, then we shall 
write $x \top \alpha$ and call the pair $(x,\alpha)$ is {\it transversal} 
(according to the context, we may also use  the terminology {\it remote}, {\it distant} or {\it in general position}). Notice that the conditions stated above simply express that $\XX$ is covered by the sets
$$
\v_\alpha := \alpha^\top :=
\{ x \in \XX  \vert \, x \top \alpha \}
$$ 
for $\alpha \in \XX'$, and
$\XX'$ is covered by the sets $\v_x'$ defined dually for $x \in \XX$.

\msk \nin {\bf 1.2. Linear and affine pair geometries.} 
Let $\K$ be a commutative unital ring.
A {\it linear pair geometry over $\K$} is a pair geometry such
that, for every transversal pair $x \top \alpha$, both $(\v_\alpha,x)$ and 
$(\v_x',\alpha)$ are equipped with a $\K$-module structure (with zero vectors $x$ and $\alpha$, respectively).
We adress the sets $\v_a$, $a \in \XX'$, as  {\it linear
parts} or {\it linear charts} of $\XX$,  and dually for $\XX'$.
We say that a linear pair geometry is {\it of commutative type}, or shorter:
is an {\it affine pair geometry}, if, for all $a \in \XX'$,
the underlying affine structure of the $\K$-modules
$(\v_a,x)$ does not depend on the choice of the element
$x \in a^\top$, and dually.
Thus we may speak of ``the'' affine space $\v_a$, and the linear
structures on $\v_a$ are related among each other by the usual formulae
$$
u +_x v = u-x+v, \quad r_x v = (1-r)x + rv
$$
where addition and scalar multiplication without subscript refers to one
fixed origin $ 0 \in \v_a$.

\msk \nin {\bf 1.3. The structure maps.}
Let $(\XX,\XX')$ be a linear pair geometry. If $(x,a), (y,a) \in (\XX \times \XX')^\top$ and $r \in \K$, then let
$r_{x,a}(y):= r y$ denote the product $r \cdot y$ in the $\K$-module
$\v_a$ with zero vector $x$. Moreover, we define the {\it structure maps} by
$$
\P_r: (\XX \times \XX' \times \XX)^\top  \to \XX, \quad
(x,a,y)   \mapsto \P_r(x,a,y) :=r_{x,a}(y),
\eqno (1.1)
$$
where
$$
 (\XX \times \XX' \times \XX)^\top  = 
\{ (x,a,y) \in \XX \times \XX' \times \XX \vert \,
x \top a, \, y \top a \},
\eqno (1.2)
$$
and dually the maps  $\P_r'$ are defined.
Similarly, we derive structure maps from the vector addition:
$$
\S: (\XX \times \XX' \times \XX \times \XX)^\top  \to \XX, \quad
(x,a,y,z)   \mapsto \S(x,a,y,z) := y +_{x,a} z,
\eqno (1.3)
$$
where the sum is taken in the $\K$-module $(\v_a,x)$ and where
the set $(\XX \times \XX' \times \XX \times \XX)^\top$ is defined
by a similar condition as in (1.2). Dually, $\S'$ is defined.

\msk \nin {\bf 1.4. Duality.}
The linear pair geometry $(\XX',\XX;\top;\P',\P,\S',\S)$ is called the
{\it dual geometry} of $(\XX,\XX';\top;\P,\P',\S,\S')$.

\msk \nin {\bf 1.5. Morphisms.\ I.}
{\it Homomorphisms} of linear pair geometries
are pairs of maps $(g,g'):(\XX,\XX') \to (\YY,\YY')$ which 
preserve transversality and which are compatible
with the multiplication and addition maps in the sense that
$g \P_r(x,a,y)=\P_r(gx,g'a,gy)$, etc.
This means simply that $g$ induces by restriction
a linear map from $(\v_a,x)$ to $(\v_{g'a},g(x))$, and
dually. In particular, we can speak of the {\it automorphism
group} $\Aut(\XX,\XX')$ of a linear pair geometry $(\XX,\XX')$.
An {\it anti-automorphism} is an isomorphism of $(\XX,\XX')$
onto its dual geometry; a {\it correlation} is an
anti-automorphism of order 2 (i.e. $g'=g^{-1}$);
a correlation is called a {\it polarity} if it admits
{\it non-isotropic points} (i.e. there is $x \in \XX$ such
that $x \top g(x)$) and a {\it null-system} if all points are isotropic.

\ssk
A {\it base point} in $(\XX,\XX')$
is a pair $(o,o') \in (\XX \times \XX')^\top$.
The {\it structure group} is the group
$\Aut(X,X';o,o')$ of automorphisms fixing the base point. 
>From the definitions it follows that this group acts linearly
on $\v_{o'} \times \v_o'$.

\msk \nin {\bf 1.6. Morphisms.\ II.}
{\it Adjoint pairs of morphisms} are given by pairs 
$g:\XX \to \YY$, $h:\YY' \to \XX'$
such that transversality is preserved in the sense that
$x \top h(a)$  iff  $g(x) \top a$, and, whenever
$(x,h(a))$ and $(y,h(a))$ are transversal, then
$$
g \P_r(x,ha,y) = \P_r(gx,a,gy)
\eqno (1.4)
$$
and similarly for $\P_r',\S$ and $\S'$; we then write $h=g^t$.
This means that $g$ induces a linear map from $(\v_{ha},x)$ to $(\v_a,gx)$. 
Note that every isomorphism $(g,g')$ gives rise to an adjoint pair
$(g,(g')^{-1})$, and conversely, every bijective adjoint pair gives
rise to an isomorphism.

\msk \nin {\bf 1.7. Connectedness.} 
We will say that two points $x,y \in \XX$ are
{\it on a common chart} if there is $a \in \XX'$,
such that $x,y \in \v_a$. 
Equivalently,  $\v_x' \cap \v_y' \not= \emptyset$. 
We will say that $x,y \in \XX$ are {\it connected}
if there is a sequence of points $x_0=x,x_1,\ldots,x_k=y$ such that
$x_i$ and $x_{i+1}$ are on a common chart.
This defines an equivalence relation on $\XX$ whose
equivalence classes are called {\it connected
components of $\XX$}. By duality, connected components
of $\XX'$ are also defined. The geometry is called
{\it connected} if both $\XX$ and $\XX'$ are connected.

\msk
\nin {\bf 1.8. Stability.}
The pair geometry $(\XX,\XX',\top)$ will be called {\it stable}
if any two points $x,y \in \XX$ are on a common chart,
and dually for any points $a,b \in \XX'$.
Clearly, a stable geometry is connected (the converse is not true).

\msk \nin {\bf 1.9. Direct products.}
The {\it direct product} of two linear pair geometries,
or of a family of linear pair geometries,
is defined in the obvious way. 

\msk \nin {\bf 1.10. Non-degeneracy.} 
The geometry is called {\it non-degenerate} if, for $a,b \in \XX'$,
 $\v_a = \v_b$ implies $a=b$, and dually.  

\msk \nin {\bf 1.11. Flat geometries.}
For any $\K$-module $V$ we may consider the geometry
$(\XX,\XX')=(V,p)$, where $\XX'=\{ p \}$ is reduced to a point,
and the multiplication map $\P_r(x,p,y)=(1-r)x+ry$
reflects just the usual affine structure of $V$.
In a similar way, we may define the flat geometry
$(\XX,\XX')=(V,V')$ associated to a pair of $\K$-modules
by letting the multiplication maps be independent
of the middle argument. These geometries are the most degenerate cases.

\sectionheadline
{2. Intrinsic subspaces in linear pair geometries}

\nin {\bf 2.1. Subspaces.}
Assume $(\XX,\XX')$ is a linear pair geometry over $\K$.
A pair $(\YY,\YY')$ of subsets
 $\YY \subset \XX$, $\YY' \subset \XX'$ is called a
{\it subspace of $(\XX,\XX')$}, if:

\msk
\item{(1)} For all $x \in \YY$ there exists an element  $a \in \YY'$ such that $x \top a$. Dually, for all $a \in \YY'$ there exists an element  $x \in \YY$ such that $x \top a$.
\item{(2)} For every $(x,a) \in (\YY \times \YY')^\top$ the set
$\YY_a:=\YY \cap \v_a$ is a {\it linear} subspace of $(\v_a,x)$ and
$\YY_x':=\YY' \cap \v_x'$ is a {\it linear} subspace of $(\v_x',a)$.

\msk \nin
Then clearly $(\YY,\YY')$ is a linear pair geometry in its own right.
Examples are provided by Lagrangian geometries, which are subspaces of Grassmann geometries, cp.\ Section 4.

\msk \nin {\bf 2.2. Intrinsic subspaces.}
A subset
 $\II \subset \XX$ is called an {\it intrinsic subspace (in $\XX$)}, 
if  $\II$ appears linearly with respect to {\it all possible}
linearizations or, more precisely, if
for all  $a \in \XX'$ with
$\v_a \cap \II \not= \emptyset$ and for all $x \in \v_a \cap \II$, the
set $(\v_a \cap \II,x)$ is a linear subspace of $(\v_a,x)$.
Given an intrinsic subspace $\II$, we put
$$
\II' := \bigcup_{x \in \II} \v_x'.
$$
Then $(\II,\II')$ is a subspace (in the sense of 2.1)
having the particular property that $\II'$ contains the whole linear part $\v_x'$ for every $x \in \II$.

\Lemma 2.3.
The intersection $\II:=\bigcap_{i \in J} \YY_i$ of an arbitrary family
 $(\II_i)_{i\in J}$ of intrinsic subspaces is an intrinsic subsace again.

\Proof. Let $x \in \II$ and $a \in \v_x$. If $i \in J$, then $(\II_i \cap \v_a,x)$ is a linear subspace of $(\v_a,x)$ and, hence,
 $\II \cap \v_x = \cap_{i \in J} (\II_i \cap \v_x)$
is a linear subspace of $(\v_a,x)$, too.
\qed

\msk \nin
{\bf 2.4. Principal intrinsic subspaces.}
Every singleton $\{ x\}$ (with $x \in \XX$) is an intrinsic subspace of
$\XX$. It is the smallest intrinsic subspace in $\XX$ containing $x$. 
For $S \subset \XX$, the subset
$$
\langle S \rangle:= \bigcap_{\II \, {\rm intrinsic \,subspace}, \,
\atop S  \subset  \II} \II,  
$$
called the {\it intrinsic subspace generated by} $S$, is an intrinsic subspace, thanks to the preceding lemma. Clearly, $\langle S \rangle$ is the smallest intrinsic subspace containing $S$.
For $x,y \in \XX$, we shall write
$$
x \lor y :=\langle \{ x,y \} \rangle  \subset X,
$$
called the {\it principal intrinsic subspace generated by
$x$ and $y$}.
If $x,y$ are not on a common chart (cf.\ 1.5), then $x \lor y = \{ x,y \}$,
and
if $x,y$ are on a common chart $\v_a$, then $x \lor y$ contains at least the
affine lines $\{ \P_r(x,a,y) | \, r \in \K \}$ and
$\{ \P_r(y,a,x) | \, r \in \K \}$ with respect to  
$(\v_a,x)$ and $(\v_a,y)$ joining $x$ and $y$
(which are the same sets if our geometry is an affine pair geometry).

\ssk
We write also $\II \land \JJ$ for the intersection of 
two intrinsic subspaces $\II$ and $\JJ$
and $\II \lor \JJ:=\langle \II \cup \JJ \rangle$ for the smallest intrinsic
subspace containing $\II \cup \JJ$.
Then the set of all intrinsic subspaces (resp. the set
of all intrinsic subspaces in $\XX$ containing a given point
$o \in \XX$) forms a {\it complete lattice with respect to inclusion and
 the operations $\land$ and $\lor$.} 

\msk \nin
{\bf 2.5. Trivial intrinsic subspaces.}
We say that an intrinsic subspace $\II$ is {\it trivial} if, 
for all $x \in \II$ and $a \top x$ the underlying affine structure of
the $\K$-module
$\m:=(\v_a \cap \II,x)$ depends neither on the choice of $x \in \m$ nor
on the choice of $a$ with $\m \subset \v_a$.

\msk \nin {\bf 2.6. Intrinsic subspaces given by images of adjoint pairs.}
Assume $(g,h):(\YY,\XX') \to (\XX,\YY')$ is an adjoint pair
of morphisms (cf.\ 1.6).
Then the image $\im(g)$ is an intrinsic subspace in $\XX$:
in fact, whenever $x,y \in \YY$ and $a \in \v_{gx}' \cap \v_{gy}'$, we have
$$
\P_r(g(x),a,g(y)) = g \P_r(x,h(a),y) \in \im(g),
$$
and similarly for $\S$; hence $\im(g)$ is an intrinsic subspace. 

\msk \nin
{\bf 2.7. Flags, chain conditions and rank.}
A {\it flag of intrinsic subspaces} is a chain
$$
\ldots \II_{-2} \subset \II_{-1} \subset \II_0 \subset \II_1 \subset \ldots 
$$
of intrinsic subspaces in $\XX$, and its {\it length},
possibly infinite, is the number of strict inclusions.
If a point $o \in \XX$ is fixed, a 
{\it principal flag centered at $o$} is given by a flag
$$
\matrix{
 o \lor x_{-2} & \subset & o \lor x_{-1} 
 & \subset & \ldots & \subset & o \lor x
\cr} 
$$
of principal intrinsic subspaces one of whose generators is $o$
and which shall all be included in some intrinsic subspace
of the form $o \lor x$.
The supremum of all lengths of (principal) flags of finite length is
called the {\it (principal) rank of $\XX$}, and
the supremum of all lengths of principal flags of finite length 
centered at $o$ and included in $o \lor x$ is
called the {\it (principal) rank}. 
We say that the geometry $(\XX,\XX')$ satisfies the
{\it descending (ascending) chain condition on principal 
subspaces (dcc, resp. acc)} if every
descending (ascending) flag of principal subspaces becomes stationary.

\msk \nin {\bf 2.8. Intrinsic lines.}
Elements of $\XX$ are called {\it points}.
A connected intrinsic subspace $\II$ in $\XX$
is called  an {\it intrinsic line} if it is not a singleton
and if all intrinsic subspaces properly contained in $\II$ are singletons. 
We will also use the term {\it minimal intrinsic
subspace}. 
If $x$ and $y$ are different points on a line $\II$, then $x\lor y=\II$.
A line is called {\it simple} if it is not a trivial
intrinsic subspace.
Thus affine lines in flat geometries of vector spaces 
are not considered as simple.
Note that, in general, intrinsic lines need not exist; 
they do exist if our geometry satisfies a descending chain condition.

\msk
\nin {\bf 2.9. A word of warning concerning the ``horizon.''}
For any $\alpha \in \XX'$, we define the {\it horizon (of the linear part
corresponding to $\alpha$)} to be the set
$$
H_\alpha := \XX \setminus \v_\alpha.
$$
In general $H_\alpha$ will {\it not} define an intrinsic
subspace in $\XX$, nor is it the first component of some subspace
$(\YY,\YY')$, and for this reason we do not use the term ``hyperplane at infinity''
which the reader might be used to from projective geometry.
In fact, one can show that essentially the only example where $H_\alpha$ is
an intrinsic subspace is the  example of ordinary projective geometry
over a field, $(\XX,\XX')=(\Gras_{1,n}(\K),\Gras_{n,1}(\K))$ (see below,
 Chapter 3).
On the other hand, in the case of geometries associated to Jordan pairs
(cf.\ Section 5), there exist non-trivial relations between
the description of
the horizon and the notion of rank.

\sectionheadline
{3. Flag geometries}

\nin In this and the next two chapters we discuss the most important
examples of linear pair geometries and their intrinsic subspaces.

\msk
\nin {\bf 3.1.} {\bf Transversality of flags.} Assume 
 $\m$ is a $\K$-module and $k \in \N$. A
 {\it flag of length  $k$}, or
{\it $k$-filtration of $\m$}, is a sequence 
 $\f=(\f_1,\f_2,\ldots,\f_{k-1},\f_k)$
of ascending linear subspaces
$$
\matrix{ 0 & \subset & \f_1 & \subset & \f_2 & \subset & \ldots 
& \subset & \f_{k-1} & \subset & \f_k=\m . \cr}
$$
In order to avoid trivialities we assume that all inclusions are strict.
Thus, the length of the flag is the number of inclusions. 
A flag of length two simply is a proper subspace. A flag of length three will be called a {\it short flag}.


\ssk
Two flags $\f$ and $\f'$ of the same length $k$ are called
 {\it transversal}, 
and we write $\f \top \f'$, 
if they are ``crosswise complementary" in the following sense:
$$
\forall i = 1,\ldots,k: \quad \quad \f_i \oplus \f_{k-i}' = \m.
$$
If $\K$ is a field, then every flag admits at least one transversal flag.
If $\K$ is just a ring, then this is not true in general.
Therefore, we will only consider {\it complemented flags}, i.e.
flags that admit a transversal flag.
The set of all complemented flags of length $k$ will oftenly denoted by 
${\cal F}={\cal F}_k$.
Obviously,
$
({\cal F},{\cal F};\top)
$
is a pair geometry in the sense of Section 1.1. For $k=2$, 
by definition this geometry is called the
{\it Grassmann-geometry} $(\Gras(\m),\Gras(\m);\top)$.

\msk \nin {\bf 3.2. Gradings and filtrations.}
A {\it $\Z$-grading of $\m$ of length $k$} 
is a family $\g=(\g_1,\ldots,\g_k)$ of submodules such that
$$
\m = \g_1 \oplus \ldots \oplus \g_{k}  .
$$
Let us denote by ${\cal G}={\cal G}_k$ 
the set of all gradings of  $\m$ of length $k$.
Each grading $\g=(\g_1,\ldots,\g_k)$ defines a pair of transversal flags 
of length $k$:
for $j=1,\ldots,k$ let
$$
\f^+_j  := \g_1 \oplus \ldots \oplus \g_j, \quad \quad
\f^-_{k-j}  := \g_{j+1} \oplus \ldots \oplus \g_k.
$$
Obviously,
$\f_j^+ \oplus \f_{k-j}^- = \g_1 \oplus \ldots \oplus \g_{k} = \m$,
whence the following flags are transversal
$$
\eqalign{
\f^+:=\f^+(\g): & \quad \quad 0 \subset \f_1^+ \subset \f_2^+  \subset \ldots \subset
\m=\f^+_k   \cr
\f^-:=\f^-(\g): & \quad \quad 0 \subset  \f_1^- \subset \f_2^- \subset \ldots \subset
\m=\f_k^-.  \cr}
$$
Observe that the grading $\g$ can be recovered from
$(\f^+,\f^-)$ via $\g_j = \f_j^+ \cap \f_{k+1-j}^-$ (for $j=1,\ldots,k$).

\Proposition 3.3.
Two flags $\e,\f$ of length $k$ in a $\K$-module $\m$
 are transversal if and only if they are derived from a grading of $\m$, i.e.,
 if and only if there exists a grading $\g$ such that
 $\e=\f^+(\g)$ and $\f=\f^-(\g)$.
If $\e$ and $\f$ are transversal, then the corresponding 
grading of $\m$ is given by
$\g_j = \e_j  \cap \f_{k+1-j}$.

\Proof.
As already remarked, it is clear that flags derived from a grading are
transversal.
Let us prove the converse by induction on $k$. 
For $k=2$ the claim is trivial.
We assume now that, for transversal flags of length $k$ in arbitrary
$\K$-modules, it has already been shown that
$\g_j = \e_j \cap \f_{k+1-j}$ defines a grading.
Given two transversal flags 
$\e$ and $\f$ of length $k+1$ in a $\K$-module $\m$, 
let us take intersections with $\n:=\e_{k}$. Putting $\f_j' := \f_j \cap \e_{k}$, we consider two flags of lenght $k$ in the space $\n$, namely
$$
\eqalign{
& \e_1 \subset \e_2 \subset \ldots \subset \e_{k-1} \subset \e_k=\n,         \cr
& 0=\f_1' \subset \f_2' \subset \ldots \subset \f_{k}' \subset \n,\cr}
$$
Let us show that these two flags are transversal, i.e., that
$$
\forall j=1,\ldots,k-1: \quad   \quad
\e_{k} = \e_j \oplus (\f_{k+1-j} \cap \e_{k}).
$$
First of all, the sum on the right hand side is direct since
$\e_j \cap \f_{k+1-j}=0$.
The inclusion ``$\supset$" is clear since $\e_j \subset \e_{k}$, and
the inclusion ``$\subset$" follows since every element $x \in \m$
can be decomposed as $x=x_1+x_2$ with $x_1 \in \e_j$ and
$x_2 \in \f_{k+1-j}$, therefore, if $x \in \e_{k}$,
it follows that $x_2 = x - x_1$ belongs both to $\f_{k+1-j}$ and
to $\e_{k}$.
Applying the induction hypothesis to these two transversal flags of 
length $k$, it follows that
$$
\e_{k} = \bigoplus_{j=1}^{k} (\e_j \cap \f_{k+2-j} \cap \e_{k}) =
\bigoplus_{j=1}^{k} (\e_j \cap \f_{k+2-j}).
$$
Finally, we have that $\m = \f_1 \oplus \e_{k} =
\oplus_{j=1}^{k+1} (\e_j \cap \f_{k+2-j})$.
\qed

\msk \nin {\bf 3.4. Groups and matrix notation.}
Given a flag $\f$ of length $k$, we define the stabilizer group $P(\f)$ (``parabolic") and its ``unipotent radical" $U(\f)$ by
$$
\eqalign{
P(\f) & = \{ g \in \Gl(\m) \vert \, \forall i : \, g(\f_i) = \f_i \},
\cr
U(\f) & = \{ g \in \Gl(\m) \vert \, \forall i : \, (g-\id)(\f_i) 
\subset \f_{i-1} \}.
\cr}
$$
For a given grading $\m=\g_1 \oplus \ldots \oplus
\g_k$, we write endomorphisms $\alpha$ of $\m$ as $k \times k$-matrices
$$
(A_{ij})_{i,j=1,\ldots,k}, \quad 
\quad A_{ij} = \pr_i \circ \alpha \circ \iota_j:
\g_j \to \g_i,
$$
with projections $\pr_i:\m \to \g_i$ and injections $\iota_j:\g_j \to \m$.
Then, if $\f$ is a fixed filtration of length $k$, 
and fixing an arbitrary transversal filtration $\e$ (hence a grading,
according to Prop.\ 3.3),
elements of the group $P(\f)$
are represented by upper triangular matrices with invertible diagonal
entries, and elements of $U(\f)$ are represented by upper triangular
matrices with diagonal entries $\id_{\g_i}$, $i=1,\ldots,k$.
(Sometimes it will also be useful to consider the flag $\e$ as ``ascending''
and $\f$ as ``descending''; then the groups $U(\f)$ and $P(\f)$ will be
represented by lower triangular matrices.)
 We define
also the ``Lie algebra of $P(\f)$, resp.\ $U(\f)$" by
$$
\eqalign{
\p(\f) & : = \{ X \in \End(\m) \vert \, \forall i : \, X(\f_i) 
\subset \f_{i} \}, \cr
\u(\f) & : = \{ X \in \End(\m) \vert \, \forall i : \, X(\f_i) 
\subset \f_{i-1} \}.   \cr}
$$
Linear operators $X \in \u(\f)$ are nilpotent: $X^k=0$. Hence,
if the integers $2,\ldots,k-1$ are invertible in $\K$, the exponential map
$$
\exp:\u(\f) \to U(\f), \quad X \mapsto \exp(X) = \sum_{i=0}^{k-1} {X^i \over
i!}
$$
is a bijection whose inverse is
given by $\log(\1 + Y)=- \sum_{i=1}^{k-1} {(-Y)^i \over i}$.

\Theorem 3.5 (The linear pair structure on flag geometries).
Let ${\cal F}_k$ be the flag geometry of all complemented flags of
length $k$ in a $\K$-module $\m$.
\item{(1)}
Given a flag $\f \in \FF_k$, the group $U(\f)$
acts simply transitively on the set $\f^\top$ of flags that are
transversal to $\f$.
\item{(2)}
Assume that $2,3,\ldots,k-1$ are invertible in $\K$. Then
the flag geometry $(\XX,\XX')=({\cal F}_k,{\cal F}_k)$ 
is a linear pair geometry over $\K$, where for a transversal pair 
$(\f,\e)$ the linear structure on $(\f^\top,\e)$ is  
defined by declaring the bijection
$$
\u(\f) \to \f^\top, \quad X \mapsto \exp(X) . \e
$$
to be a linear isomorphism. 
\item{(3)} 
The linear pair geometry defined in part (2) 
is an affine pair geometry if 
$(\XX,\XX')$ is a Grassmann geometry, i.e.\ if $k=2$.

\Proof. (1)
We prove the claim by induction on the length $k$ of the flag:
for $k=2$ (complements of a single subspace) the assertion
is well-known: in this case $\f^\top$ is an affine space on which the
vector group $U(\f)$ acts as translation group.
Next, suppose that the claim is proved for flags of length $k$. Consider 
a flag $\f$ of length $k+1$ together with two flags
$\e$ and $\e'$ which are transversal to $\f$.
As both $\e_1$ and $\e_1'$ are complements of $\f_k$,
we find $u_1 \in \Gl(\m)$ with $u_1(\e_1)=\e_1'$ and
$u_1\vert_{\f_k} = \id_{\f_k}$. 
In particular, $u_1$ is an element of $U(\f)$. 
Now consider the flag $\tilde \f$ of length $k$
$$
0 \subset \f_1 \subset \ldots \subset \f_k =:\n.
$$
Both flags $\tilde \e_j:=\e_j \cap \f_k$  and
$\tilde \e_j':=\e_j' \cap \f_k$ ($j=2,\ldots,k+1$)
 are transversal to $\tilde \f$, cp.\ see proof of Prop.\ 3.3.
By the induction hypothesis, there exists an element
$u_2 \in U(\tilde f)$ with $u_2 (\tilde \e) = \tilde \e'$.
Therefore, $u_3 := \id_{\e_1} \oplus u_2$ belongs to $U(\f)$ and satisfies 
$u_3(\tilde \e_j) = \tilde \e_j'$.
In matrix form, the elements $u_1$ and $u_3$ of $U(\f)$ are given by
$$
u_1 = \pmatrix{\id_{\e_1} & ** \cr 0 & \id_{\f_k} \cr}, \quad \quad
u_3 = \pmatrix{\id_{\e_1} & 0 \cr 0 & u_2 \cr}.
$$
Since $\e_j = \tilde \e_j \oplus \e_1$ and
 $\e_j' = \tilde \e_j' \oplus \e_1$, 
we have
$u_3(\tilde \e_j) = \tilde \e_j'$ for $j=2,\ldots,k$.
Finally, $u := u_1 \circ u_3$ is an element of 
$U(\f)$ such that $u(\e)=\e'$.

\ssk
Uniqueness:
Assume $u \in U(\f)$ is such that $u(\e)=\e$.
Then the matrix of $u$ clearly is the identity matrix, and hence
$u=\id_\m$.

\msk (2)
The exponential map $\u(\f) \to U(\f)$ being a bijection,
composition with the bijection $U(\f) \to \f^\top$, $g \mapsto g.\e$
from part (1) yields a bijection between $\f^\top$ and the $\K$-module
$\u(\f)$, which we use to define a linear structure on $(\f^\top,\e)$.
Similarly, we define the linear structure on $(\e^\top,\f)$,
and thus define a linear pair geometry.

\ssk (3)
The linear pair geometry just defined 
is an affine pair geometry if 
the groups $U(\f)$ are all abelian, and this is the case if $k=2$.
\qed

\msk \nin {\bf 3.6. Fixing the type.}
In general, the whole flag geometry $(\FF,\FF)$ is neither connected
(cf.\ 1.7) nor homegeoneous. Similar to Grasmmanian geometry, we therefore 
fix a {\it type} $\o$ (some flag of length $k$) and a {\it co-type}
$\o'$ (some transversal flag), and then put 
$$
\eqalign{
\XX & := \FF_\o^{\o'} := \{ \f \in \FF \vert \, \f \cong \o, \, \exists \f': 
\f' \cong \o', \f' \top \f \}, \cr
\XX' & := \FF_{\o'}^\o :=
\{ \f \in \FF \vert \, \f \cong \o', \, \exists \f': 
\f' \cong \o, \f' \top \f \},   \cr}
$$
where $\f \cong \e$ means that $\f_j$ and $\e_j$ are isomorphic as
$\K$-modules for all $j$. 
Note that, if $\f$ is some flag, then all transversal flags $\f'$ of
$\f$ have necessarily  the
same type and co-type. Thus, if the intersection of $\f^\top$ with $ \FF_{\o'}^\o $ is non-empty, then necessarily the whole set
$\f^\top$ belongs to $ \FF_{\o'}^\o $ and, consequently,
$(\XX,\XX')=( \FF_\o^{\o'}, \FF_{\o'}^\o) $ forms a subgeometry of $(\FF,\FF)$.
If $\K$ is a field and $\m = \K^n$, then
$\XX$ and $\XX'$ are spaces of flags of a given sequence of dimensions,
and if moreover $k=2$ we get the usual Grassmannians of $p$-, resp.\
$q$-spaces in $\K^n$.

\msk \nin {\bf 3.7. Remark: the structure of the structure maps.}
In the case of the Grassmann geometry one can give a fairly explicit
formula for the structure maps $\P_r$, and one can find algebraic 
identities that, in some sense, characterize these maps (see [Be02], [Be04]).
For $k>2$, no results of this kind seem to be known.
In this case, the linear parts
$\alpha^\top$ carry many different affine structures (since the linear spaces
$(\alpha^\top,x)$ for various choices of the origin $x$ in $\alpha^\top$ 
have no longer the same underlying affine space).  
A situation where one has to take account of many affine structures on
the same underlying set
occurs also in the theory of linear connections as developed in [Be06] --
we believe that these similarities are no fortuities.
The collection of the various affine structures on $\alpha^\top$ contains
more information than just the group structure on $\alpha^\top$, 
and one may expect that, as in the case $k=2$, this richer structure
finally allows to develop the corresponding 
theory over rings and in infinite dimension,
where the combinatorial or Lie theoretic information alone is too weak.

\Theorem 3.8 (Standard intrinsic subspaces of flag geometries).
Assume $(\XX^+,\XX^-)$ is the flag geometry $(\FF,\FF)$ of complemented flags
of length $k$ in $\m$ 
(or one of its subgeometries obtained by fixing type and
co-type of the flags). 
We fix an index $j \in \{ 0,\ldots,k-1 \}$ 
and some subspace $\e \subset \m$. Then the set
$$
\II_{\e;j}  := \{ \f \in \XX^+ \vert \, \f_j \subset \e  \subset \f_{j+1} \}
$$
is an intrinsic subspace of $\XX^+$. Intrinsic subspaces of this type, as well as the intersections of finitely many of them, will be
called {\rm standard intrinsic subspaces} of the flag geometry.

\Proof.                                                
Fixing terminology, we say that a flag $\f$ is {\it governed by the subspace
$\e$ in position $j$}, written $ \f {\buildrel j \over \prec} \e$,
if $ \f_j \subset \e  \subset \f_{j+1} $.
Let us prove that the set $\II :=\II_{\e,j}= \{ \f \in \X | \, \f
 {\buildrel j \over \prec} \e \}$ of flags 
 governed by the subspace $\e$ in position $j$
 is an intrinsic subspace. 
If $\II$ is empty, there is nothing to prove.
Thus, let $\o \in \II$ and pick up some
$\o' \in \XX^-$ with $\o \top \o'$.
Then the theorem follows from the following

\msk
\nin {\bf Claim 3.9.} {\it
The set
$$
I:= \{ u \in U(\o') | \, u.\o {\buildrel j \over \prec} \e \}
$$
is a subgroup of $U(\o')$; more precisely, this group is the intersection
of $U(\o')$ with the ``parabolic subgroup'' $P(\e)$ stabilizing the subspace
$\e$.}

\msk \nin
Indeed, if Claim 3.9 is true, then the simple transitivity of the group $U(\o')$ on $\v_{\o'}$ shows that its subgroup $I$ acts simply
transitively on  $\II \cap \v_{\o'} = I.\o$.
As the ``Lie algebra'' of $I$ (i.e., the image of $I$ under the logarithm
$\log: U(\o') \to \u(\o')$)  is a linear space, given by 
$\u(\o') \cap \p(\e)$, the set $I.\o$ is a linear subspace of the linear space defined in Theorem 3.5 (2).

\ssk It remains to prove Claim 3.9. Notice that the relation
$ \o {\buildrel j \over \prec} \e$ is invariant under elements $g$  of 
the projective group of $\m$, i.e.\ $ \o {\buildrel j \over \prec} \e$ holds if and only if $ g.\o {\buildrel j \over \prec} g.\e$. Therefore, $u\in P(\e)$ implies that $ u.\o {\buildrel j \over \prec} \e$ and we derive that $P(\e) \cap U(\o') \subset I$. We shall show the converse inclusion 
$I \subset P(\e) \cap U(\o')$ by an explicit ``matrix calculation'': let $u$ be an element of $I$. We have to prove that $u$ is an element of $P(\e)$, i.e.\ that $u(\e) \subset \e$. As $u$ is an element of $U(\o')$, we may represent $u$ by a {\it lower} triangular matrix\footnote{$^1$}{\eightrm
Thus, thinking in terms of matrices, we consider $\o$ as an ascending flag,
and $\o'$ as a descending flag; since $U(\o')$ preserves the descending flag,
it is represented by lower triangular matrices. The flag $0 \subset \e
\subset \m$ should also be considered as ascending.}
$(A_{i,j})$ with respect to the grading 
$\m = \g_1 \oplus \ldots \oplus \g_k$ corresponding to $(\o,\o')$.
As in the proof of Prop.\ 3.3, we have
that $\e = \o_j \oplus (\e \cap \g_{j+1})$. Recall that the condition $u.\o {\buildrel j \over \prec} \e$ means that

\ssk
\item{a)}
$u.\o_j =(u.\o)_j \subset \e$, and
\item{b)}
$\e \subset (u.\o)_{j+1} = u.  \o_{j+1}$.

\ssk \nin In terms of matrices, Condition a) can be rephrased as follows:
we represent $u$ by a $3 \times 3$-matrix with respect to the decomposition 
$\m = (\g_1 \oplus \ldots \oplus \g_j) \oplus \g_{j+1} \oplus
(\g_{j+2} \oplus \ldots \oplus \g_k)$.  From the assumption $\o_j \subset \e \subset \o_{j+1}$ (i.e.\ $(\g_1 \oplus \ldots \oplus \g_j) \subset \e
\subset (\g_1 \oplus \ldots \oplus \g_{j+1})$), we see that condition a) is fulfilled if and only if $u$ is represented by a matrix
$$
u = \pmatrix{\1 & 0 & 0 \cr B & \1 & 0 \cr 0 & * & \1 \cr}, \quad \quad
B = \pmatrix{* & \ldots & * \cr 0 & \ldots & 0 \cr},
$$
where the matrix $B$ just symbolically visualizes the condition $\im (A_{j+1,s}) \subset \e$, $1 \leq s < j$.
We claim that, similarly, condition b) is equivalent to
$$
u = \pmatrix{\1 & 0 & 0 \cr * & \1 & 0 \cr 0 & C & \1 \cr}, \quad \quad
C = \pmatrix{0 & * \cr & \vdots \cr 0 & * \cr}
$$
where the matrix $C$ stands for the condition
$\g_{j+1} \cap \e \subset \ker (A_{s,(j+1)})$ for all $s \leq j$.
In fact,
if $\m$ is free over $\K$ of finite rank $n$, then the assertion follows from a) by transposing, i.e.\ by applying the preceeding 
case to the dual space (where now ascending and descending flags change their r\^ole). In the general case, we argue as 
follows\footnote{$^2$}{\eightrm
One should expect that
there is a general ``duality principle'' allowing to deal with such situations in
an axiomatic context, such as the duality principle known in Jordan theory, see
[Lo75, I Prop.\ 2.9]}: 
if $y \in (\g_{j+1} \cap \e)$, then we have to show that $y=uy$.
Since $\e \subset \u.\o_{j+1}$, 
 every element $y \in (\g_{j+1} \cap \e)$ has
the form $y=u.x=u.\sum_{i\leq j+1} x_i=\sum_{i \leq j+1} x_i + 
\sum_{i \leq j+1 \atop s>i}
A_{is} x_i$ (where $x_i =\pr_i (x) \in \g_i$) for some $x \in \o_{j+1}$. 
As $y \in \g_{j+1}$ holds by assumption, this 
triangular system has the unique solution
$x_{j+1}=y, x_1=0,\ldots,x_j=0$, whence $x=y$ and, hence, $y=uy$ as had to be
shown.

\ssk
Combining conditions a) and b), it is easy to see that $u$ satisfies the condition $u.\e \subset \e$, and hence belongs to $U(\o') \cap P(\e)$. This finishes the proof of the claim.
\qed

\ssk
\nin Typically, intersections of intrinsic subspaces of the type 
just  described lead to the
following situation: let $\e$ be a flag of length $k+1$; then the set
$$
\II = \{ \f \in \XX | \, \forall j = 1,\ldots,k: \e_j \subset \f_j \subset \e_{j+1} \}
\eqno (3.1)
$$
is an intrinsic subspace. Again, its linear image is the intersection of a group
of unipotent lower triangular matrices with the ``(upper triangular)
parabolic subgroup'' stabilizing the flag $\e$. 
                                                   
\msk \nin {\bf 3.10.  Case of a Grassmann geometry.}
Let $(\XX,\XX')=(\Gras_\o^{\o'}(\m), \Gras_{\o'}^\o(\m))$ 
be a Grassmann geometry of fixed type and co-type. Fixing a flag of length
three, $\e=(0 \subset \e_1 \subset \e_2 \subset \m)$, we consider the space
$$
\II_\e : = \II_{\e_1;0} \cap \II_{\e_2;1} =
 \{ x \in \Gras(\m) \vert \, \e_1 \subset x \subset \e_2 \}.
\eqno (3.2)
$$
With respect to a base point $(\o,\o')$, where $\o \in \II_\e$ 
(i.e.,  $\e_1 \subset \o \subset \e_2$), the space $I_\e = \II_e \cap V$
is then given as follows:
The group $U(\o')$ is the abelian vector group
$$
U(\o') = \{ u=\pmatrix{\1 & 0 \cr f & \1 \cr} \vert \, f \in \Hom(\o,\o') \}
$$
Then $u.\o = \Gamma_f$ can be identified with the graph 
$\Gamma_f = \{ \bigl({u \atop fu }\bigr) \vert \, u \in \o \}$  of $f$,
and (writing $\e_2 = \o \oplus \e_2'$ with $\e_2'=\e_2 \cap \o'$)
the compatibility conditions $\e_1 \subset \Gamma_f$, 
$\Gamma_f \subset \e_2$, are equivalent to
$$
\e_1 \subset \ker(f), \quad \im(f) \subset \e_2'.
$$
In other words, $f \in I_{\e_1,\e_2'} = \{ f \in \Hom(\o,\o') \vert \,
\e_1 \subset \ker(f), \, \im(f) \subset \e_2' \}$.
Clearly, this is a linear subspace
of $\Hom(\o',\o)$. Moreover, it is an {\it inner ideal} of $\Hom(\o',\o)$,
see Appendix A, Section A.3. We will prove that, under suitable assumptions,
all intrinsic subspaces are obtained in this way:

\Theorem 3.11.
Assume that $\K$ is an infinite field of characteristic not $2$ and 
that $(\XX^+,\XX^-)$ is the
 Grassmannian geometry 
$(\Gras_{p,q}(\K),\Gras_{q,p}(\K))$  
 of $p$- and $q$-spaces in a finite-dimen\-sion\-al 
vector space $\m \cong \K^{p+q}$. Then all intrinsic subspaces of $\XX^+$
have the form $\II_\e$ described in $(3.2)$ for some flag of length three,
$\e=(\e_1 \subset \e_2 \subset \m)$, such that
$\dim \e_1 \leq p \leq \dim \e_2$.
The principal subspaces are precisely the intrinsic subspaces $\II_\e$ with
$p - \dim \e_1 = \dim \e_2 - p$.

\Proof.
Suppose that $\II \subset \XX^+$ is an intrinsic subspace; we may assume that
it contains the base point $\o = \K^p$. Let $\o' \top \o$; we may 
assume that $\o'$ is the second component in the decomposition
$\K^n = \K^p \oplus \K^q$.
Let 
$$
(V^+,V^-):=(\v_{\o'},\v_\o') = (\Hom(\o,\o'),\Hom(\o' ,\o)) 
= (M(q,p;\K),M(p,q;\K))
$$
be the linearization of $(\XX^+,\XX^-)$ corresponding to this base point.
Then the affine part $I:=\II \cap V^+$ of $\II$ is a linear subspace
of $V^+$.

\ssk
We claim that $I$ is an {\it inner ideal} of $V^+=\Hom(\o,\o')$, i.e.,
$Q^+(I)V^- \subset I$, where $Q^+(X)Y=XYX$
(cf.\ Appendix A, Section A.3).
In order to prove this claim, 
let $X \in V^+$ and $Y \in V^-$  such that $X$ and $-X$ lie 
in $\v_Y$, i.e., $X \top Y$ and $-X \top Y$. 
(The set of such pairs $(X,Y)$ is non empty and Zariski-open in
$V^+ \times V^-$. This is seen by elementary arguments, or by using  
Jordan theory, noticing that that this condition is equivalent
to saying that the pairs
$(X,Y)$ and $(-X,Y)$ are {\it quasi-invertible}; in this case this means 
simply that 
$\1 - XY$, $\1 - YX$, $\1 + XY$ and $\1 + YX$ are invertible matrices.)
By a direct (but somewhat tricky) calculation one can show that
the geometric midpoint of $X$ and $-X$ in the affinization by $Y$
is given by the formula
$$
\P_{1 \over 2} (X,Y,-X) = XYX = Q^+(X) Y.
$$
(see [Be04, Eqn. (2.15)]). 
Since $\II$ is an intrinsic subspace, it follows that, whenever $X$ 
and $-X$ belong to $I=\II \cap V$ and 
$Y \in \XX^-$ is any affinization such that $X \top Y$ and $-X \top Y$,
the geometric midpoint $\P_{1 \over 2} (X,Y,-X)$ belongs again to $\II$.
Therefore, the preceding calculation shows that, whenever $X \in I$, then,
for some Zariski-dense set of elements $Y \in V^-$, the element
$Q^+(X)Y$ again belongs to $I$.
Since $\K$ is infinite, it follows that $Q^+(X)V^- \subset I$,
i.e., $I$ is an inner ideal.

\ssk 
Now, from Proposition A.4, we know that 
all inner ideals are of the form $I_{E,F}$ with subspaces $E \subset \o$,
$F \subset \o'$. Comparing with the description of the affine picture
of the intrinsic subspace $\II_\e$ above, we see that, for the flag
$\e=(E \subset \o \oplus F \subset \m)$, the intrinsic subspaces
$\II_e$ and $\II$ have the same affine picture, i.e., the same intersection
with $V^+$. 
We claim that $\II_\e$ and $\II$ also have the same ``points at infinity",
i.e., that we have in fact equality $\II = \II_\e$.
In fact, the preceding arguments show that, with respect to {\it any} 
affinization, $\II$ agrees with some standard intrinsic subspace.
Thus, covering $\II$ by a finite number $n$ of affine charts, 
the affine pictures coincide with those of spaces $\II_{\e_i}$
defined by $n$ short flags
$\e^1,\ldots,\e^n$. Since the intersection of two chart domains is
Zariski-dense in each of the affine parts and the flag is uniquely
determined by the affine picture on this intersection, it follows then
that $\e^i = \e^j=\e$ for all $i,j$, and hence $\II = \II_\e$.
\qed

\ssk In particular, the theorem says that, under suitable assumptions,
the collection of all intrinsic subspaces in a Grassmann geometry forms
again a linear pair geometry, namely a flag geometry (of short flags).
We do not know whether an analogue of  this theorem holds for $k>2$, the problem
being that we do not have enough
 knowledge on the structure maps of flag geometries that are
not ``of commutative Jordan type'' -- cf.\ Section 3.7.

\sectionheadline
{4. Lagrangian flag geometries}

\nin {\bf 4.1. Lagrangian flags.} 
 Assume $\beta$ is a strongly non-degenerate quadratic or skew-symmetric form
on $\m$, i.e., the map $\m \to \m^*$, $v \mapsto \beta(v,\cdot)$
is not  only injective, but also surjective. 
A {\it Lagrangian flag (of length $k$)} is a flag $\f$ such that
$$
\forall j=1,\ldots,k: \quad \quad \f_j^\perp = \f_{k-j}
$$
We say that a Lagrangian flag $\f$ is {\it complemented} if there exists
a transversal Lagrangian flag $\e$, and we denote by $\LL_k = \LL_k(\m)$ the
set of all complemented Lagrangian flags of length $k$ in $\m$.

\ssk Let us first discuss the 
case  $k=2$ which is particularly important.
In this case,  $\f = (0 \subset \f_1 \subset \m)$, and $\f_1$
is a {\it Lagrangian subspace}, i.e., $\f_1 = (\f_1)^\perp$, which admits
a complementary Lagrangian subspace $\e_1$. 
We call $\LL_2$ also the 
{\it Lagrangian Grassmannian} (manifold of all Lagrangian subspaces).
In matrix form, $\beta$
can then be represented on $\m=\f_1 \oplus \e_1$ by
$$
\pmatrix{ 0 & \alpha \cr \gamma & 0 \cr},
$$
with pairings $\alpha:\e_1 \times \f_1 \to \K$,
$\gamma:\f_1 \times \e_1 \to \K$.
Then  two cases have to be distinguished:

\ssk
\item{(sy)} $\beta$ is symmetric, i.e.\ $\alpha = \gamma$ 
(under identification of
$\e_1 \times \f_1$ with $\f_1 \times \e_1$)
\item{(sw)} $\beta$ is skew-symmetric, i.e. $\alpha=-\gamma$.

\ssk \nin
If $\K$ is a field and $\m = \K^n$, then necessarily $n$ is even,
$n=2m$, and for $\alpha$ we may choose the standard scalar product
on $\K^m$. The case
(sym) is then also called the {\it artinian case} and the case
(skew) the {\it symplectic case}.           
For general Lagrangian flags of length $k \in \N$ we distinguish :

\ssk
\item{(1)} $k$ is even. Then $\f_{k \over 2}^\perp = \f_{k \over 2}$
is again a Lagrangian subspace, and
all $\f_j$ with $j \leq {k \over 2}$ are isotropic subspaces.
Lagrangian flags can be constructred by choosing 
a flag of $k\over 2$ isotropic subspaces, and completing this 
flag by adding the $k \over 2$ orthogonal complements.
Note that, if $\K$ is a field, $\m = \K^n = \K^{2m}$, and
(as in case $k=2$) $\beta$ will be artinian or symplectic.
\item{(2)}
$k$ is odd. In this case
it is not necessary that Lagrangian subspaces exist; therefore the form
may be of general signature (e.g., for a Lorentzian form: 
a Lagrangian flag must be of the form $\K v \subset v^\perp \subset \m$
for an isotropic vector $v$). 


\Theorem 4.2.
The map
$$
(\FF_k,\FF_k) \to (\FF_k,\FF_k), \quad
(\e,\f) \mapsto (\e^\perp,\f^\perp)
$$
is an automorphism of the flag geometry $(\FF_k,\FF_k)$ whose fixed
point set is $(\LL_k,\LL_k)$. In particular,
the Lagrangian geometry $(\LL_k,\LL_k)$ is a subspace of the linear pair 
geometry $(\FF_k,\FF_k)$.
The linear pair geometry $(\LL_k,\LL_k)$ is an affine pair geometry if $k=2$.

\Proof.
By assumption, the map $b:\m \to \m^*$, $v \mapsto v^* = \beta(\cdot,v)$
is bijective. Since a decomposition $\m = \a \oplus \b$ gives rise to a
decomposition $\m^* = \a^* \oplus \b^*$ (with $\b^*$ corresponding to
$\a^\perp$ under the bijection $b$), it follows that 
$\e \top \f$ if and only if $\e^\perp \top \f^\perp$.
Moreover, since $b$ is a bijection, the canonical map
$\m \to (\m^*)^*$ into the double dual is also an isomorphism, and from this it
follows that $(f^\perp)^\perp = \f$.

Hence our map is an involution of the underlying pair geometry.
It is an automorphism of linear pair geometries:
for all $g \in \Gl(\m)$ and all flags $\e$, we have
 $(g.\e)^\perp = (g^*)^{-1} \e^\perp$ where $g^*$ is the adjoint of $g$
(the dual map of $g$ under the identification $b$). 
Thus
the map
$$
V^+ = \exp(\u(\o')).\o \to (V^+)^\perp = \exp(\u((\o')^\perp).\o^\perp, 
\quad u.\o  \mapsto (u^*)^{-1}.\o^\perp
$$
is linear under the exponential map since it is given by
$X \mapsto - X^*$.


Finally, since the fixed space of an automorphism is
a subgeometry, it follows that $(\LL_k,\LL_k)$ is a subgeometry.
For $k=2$ it is an affine pair geometry since it is a subspace of an
affine pair geometry.
\qed

\msk \nin {\bf 4.3. Symmetric matrices.}
Let us now focus on the case $k=2$.
With respect to a fixed base point $(\o,\o')$ (i.e., a transversal
 pair of Lagrangian subspaces), the abelian group
$$
V^+ = \UU(\o') \cap \OO(\beta) =
\{ \pmatrix{ 1 & v \cr 0 & 1 \cr} | \, v \in \Hom(\o,\o) \} \cap \OO(\beta)
$$
is described as follows:
using that $\beta$ defines a bijection $b:\m \to \m^*$,
 we get two bijections of $\o$ with  $\o'$, one of which is considered as
an idenfication, i.e, we write $\m = \o \oplus \o$. Then,

\ssk
\item{(sy)} if $\beta$ is symmetric, both bijections coincide, i.e., we write
$\beta =  \pmatrix{ 0 & \alpha \cr \alpha & 0 \cr}$; then
$V^+ \cong \Asym(\o,\alpha) $ is given by skewsymmetic operators $v$.
\item{(sw)} if $\beta$ is skew-symmetric, we write 
$\beta = \pmatrix{ 0 & -\alpha \cr \alpha & 0 \cr}$; then
$V^+ \cong \Sym(\o,\alpha) $ is given by symmetric operators $v$.
\ssk

\nin
In particular, symplectic Lagrangian geometries correspond to usual symmetric
matrices, and artinian Lagrangian geometries correspond to usual skewsymmetric
matrices.
In the case $k>2$, similar (but more complicated) descriptions of the linear model
space of the Lagrangian flags geometries can be given.

\msk \nin {\bf 4.4. Standard intrinsic subspaces.}
Intersection of the standard intrinsic subspaces of the flag geometry
$(\FF,\FF)$ from Theorem 3.11
with the subgeometry $(\LL,\LL)$ defines intrinsic subspaces in $\LL$.
For instance, for a fixed flag $\e$ of length $k+1$, intersection of the intrinsic
subspace defined by Equation (3.1) show that
$$
\LL_\e = \{ \f \in \LL \vert \, \forall j=1,\ldots,k: \,
\e_j \subset \f_j \subset \e_{j+1} \} 
\eqno (4.1)
$$
is an intrinsic subspace of $\LL$. However,
 a condition
$\e_j \subset \f_j$ automatically yields $\f_{k-j}  = \f_j^\perp \subset
\e_j^\perp$, and therefore the intrinsic subspaces $\LL_\e$ are uninteresting
unless $\e$ itself was already a Lagrangian flag.
For simplicity, assume now that $k=2$ and
 let $\e = (0 \subset \e_1 \subset \e_1^\perp \subset \m)$ be a short
Lagrangian flag. Then the affine image of the intrinsic subspace
$$
\LL_\e = \{ \f = (0 \subset \f_1 \subset \m) \in \LL \vert \, \,
\e_1 \subset \f_1 \subset \e_1^\perp \} 
\eqno (4.2)
$$
is given as follows : write $\o = \e_1 \oplus \a$; then
$$
V_\e = \{ f \in \Sym(\o) | \, \a \subset \ker f \} =
\{ f \in \Sym | \, \im f \subset \a^\perp \}
$$
(last equality: since $f$ is symmetric, we have
$\im(f) = \ker(f)^\perp$, and hence the two conditions are equivalent).
Then $V_e$ is an inner ideal in the Jordan algebra $\Sym(\o)$
(see Appendix A, A.12).
If $\K$ is a field and $\o = \K^n$, then one can show that every inner ideal
is of this form (Appendix A.12), 
and we get the following analog of Theorem 3.11:

\Theorem 4.5.
Assume that $\K$ is an infinite field of characteristic not $2$ and 
that $(\XX^+,\XX^-)=(\Lag_2,\Lag_2)$
is a finite-dimensional symplectic Lagrangian geometry.
  Then all intrinsic subspaces of $\XX^+$
are of the form $\II_\e$ given by $(4.2)$, for some Lagrangian flag of length three,
$\e=(\e_1 \subset \e_1^\perp \subset \m)$, such that
$\dim \e_1 \leq n$.
In this case, all intrinsic subspaces are principal.
\qed

Under the assumptions of the theorem,
the geometry of intrinsic subspaces is essentially isomorphic to the 
geometry $(\LL_3,\LL_3)$ (subject to dimension conditions) and hence
can be equipped with the structure of a linear pair geometry.
Its connected compenents are homogeneous under the symplectic group
$\Sp(n,\K)$, with stabilizer groups being parabolic subgroups
stabilizing some Lagrangian flag of length $3$.
Finally, a similar resulat in the artinian case holds.

\sectionheadline{5. Geometries of Lie algebra filtrations}

\nin {\bf 5.1. Filtered and graded algebras.}
Assume that $\m$ carries the structure of an {\it algebra over $\K$},
with bilinear product denoted by $x \cdot y$.
(It may be any product; only later we will assume that it is a Lie bracket
and then we will write $[x,y]$.)
Then we consider flags that are {\it filtrations} of this algebra,
i.e., they are compatible with the product in the sense that
$\f_i \cdot \f_j \subset f_{i+j}$.
Similarly, we consider only gradings $\m = \oplus_i \g_i$ which are
compatible with the algebra structure:
$\g_i \cdot \g_j \subset \g_{i+j}$.
Note that, at this point,
the choice of our indix set becomes important: so far we needed only
that the index set is totally ordered; now we need
 a structure of abelian semigroup on the index set. 
Therefore we will no longer use $\{ 1, \ldots, k \}$ as index set,
but rather speak of {\it finite $\Z$-gradings} 
$\m = \oplus_{i \in \Z} \,  \g_i$,
i.e., $\g_i \not= 0$ only for finitely
many indices, and similarly {\it finite $\Z$-filtrations} are
defined.

\Proposition 5.2. Let $\m$ be an algebra over $\K$. Then the following
objects are in canonical bijection with each other:
\item{(1)} pairs of transversal finite $\Z$-filtrations of $\m$
\item{(2)} finite $\Z$-gradings of $\m$
\item{(3)}  derivations $D:\m \to \m$ which are diagonalizable
over $\Z$ and have only finitely many eigenvalues which are contained
in some intervall $-n, \ldots ,n$ such that $1,\ldots,3n \in \K^\times$

\Proof.
It is clear that gradings give rise to pairs of transversal 
filtrations, and  the converse follows from Prop.\ 3.3:
with the notation introduced there,
two transversal filtrations $\f,\e$ give rise, first of all 
to a $\K$-module grading of $\m$ (Prop.\ 3.3), which is then easily 
seen to be an algebra grading.

\ssk
Given a grading $\m = \oplus_{i \in \Z} \, \g_i$, one checks that
the linear map $D$ given by $D(x)=ix$ for $x \in \g_i$ is a derivation
of the algebra $\m$. Conversely, given a derivation as in (3),
one checks that
the eigenspace decomposition of $D$ is a finite $\Z$-grading.
(The condition on the eigenvalues is to ensure that the eigenspace 
decomposition is direct.)
\qed

\nin
In general, the geometry of finite $\Z$-gradings of an algebra
will not carry the structure of a linear pair geometry.
However, under suitable assumptions this is the case.
Let us assume that $\m$ is a Lie algebra with bracket
$[x,y]$, and let us consider
{\it $3$-filtrations}:  these
are filtrations of length 3, with index set $\{ 1,0,-1 \}$,
i.e., filtrations of the form
$$
\eqalign{
& 0 \subset \f_1 \subset \f_0 \subset \f_{-1} = \m, \cr
& [\f_1 , \f_1] =0, \,
[\f_0 , \f_0] \subset \f_0, \,
[\f_1 , \f_0] \subset \f_1, \,
[\m , \f_1 ]\subset \f_0,\,
[\f_1 ,\m ] \subset \f_0.  \cr}
$$
Thus $\f_0$ is a subalgebra, and $\f_1$ is an abelian ideal in $\f_0$; 
moreover, left multiplication $\ad(x)$ by elements $x \in \f_1$
is three-step nilpotent. Thus, assuming that $2$ and $3$ are invertible
in $\K$, there is a well-defined automorphism $e^{\ad(x)}$ of $\g$.
In [BeNe04] the following is proved:

\Theorem 5.3.
Assume that $\m$ is a Lie algebra and $2$ and $3$ are invertible in $\K$.
Then the geometry of all inner $3$-filtrations of $\g$ carries
a natural structure of an affine pair geometry over $\K$. More precisely,
for every inner $3$-filtration $\f=(\f_1 \subset \f_0 \subset \f_{-1}=\m)$,
the vector group $\f_1$ acts simply transitively on $\f^\top$ via
$(x,\e) \mapsto e^{\ad(x)} (\e)$.
\qed

\nin Here, the term ``inner" means the following: 
a transversal pair $(\e,\f)$ of 3-filtrations is called {\it inner},
if the corresponding 3-grading $\m=\g_1 \oplus \g_0 \oplus \g_{-1}$
is of inner type, i.e.\ the derivation $D$ corresponding to this
grading via Prop.\ 5.2 is an inner derivation: 
$D=\ad(E)$ for some element $E \in \m$  (which is then called an
{\it Euler operator}).
Without this assumption the result does not hold: take, e.g., an abelian
Lie algebra; then filtrations are just ordinary flags of length three;
but the geometry of flags of length three is not an {\it affine} pair
geometry.
--
We conjecture that also for general inner filtrations of the
form $0 \subset \f_k \subset \f_{k-1} \subset \ldots \subset
\f_{-k}=\m$ the (unipotent) group $\exp(\f_1)$ acts simply transitively
on $\f^\top$, 
and that, therefore, such filtrations give rise to linear
(non affine in general) pair geometries. 
This would be very interesting since a good deal of the theory
from [BeNe04] and [BeNe05] 
would then also generalize, thus defining manifold structures on a class of 
flag geometries in arbitrary dimension and over general base fields
and rings.

\msk \nin {\bf 5.4. Fixing the type.}
As in Section 3.6, we restrict our attention to subgeometries
of the geometry of all inner $3$-filtrations by fixing a type.
In this case, this amounts to fix a base point $(\o^+,\o^-)$,
which is nothing but a fixed pair of transversal inner $3$-filtrations,
or, what is the same by Prop.\ 5.2, a fixed inner $3$-grading
$\m=\g_1 \oplus \g_0 \oplus \g_{-1}$, or a fixed Euler operator $E$.
Then the subgroup $\PE(\m,E)$ of $\Aut(\m)$ 
generated by $U^+=\exp(\g_1)$ and by $U^-=\exp(\g_{-1})$ is called the
{\it elementary projective group of $(\m,E)$}. The pair of orbits
$$
(\XX^+,\XX^-), \quad \quad \XX^\pm := \PE(\m,E) . \o^\pm
$$
is then a subgeometry of the affine pair geometry from Theorem 5.3,
called the {\it projective geometry of the $3$-graded Lie algebra
$(\m,E)$}. 

\msk \nin {\bf 5.5. Friendly Jordan pairs.}
If $\m=\g_1 \oplus \g_0 \oplus \g_{-1}$ is a $3$-graded Lie algebra
(over a ring $\K$ in which $2$ and $3$ are invertible),
then the pair $(V^+,V^-)=(\g_1,\g_{-1})$ with the trilinear maps
$T^\pm(x,y,z):=- [[x,y],z]$
is called the {\it associated Jordan pair} (see Appendix B, Section B.2).
It satisfies the identities of an abstract {\it linear Jordan pair}.
Conversely, from a linear Jordan pair one can reconstruct
a $3$-graded Lie algebra.
We say that a Jordan par is {\it friendly} if it has the following property:
for all $a \in V^\mp$, the sets
$$
U_a:=\{ x \in V^\pm \vert \, B^\pm(x,a) \in \Gl(V^\pm) \}
$$
(where $B^\pm(x,a) \in \End(V^\pm)$ is the Bergmann-Operator, see B.2)
generate $V^\pm$ as a $\K$-module; an inner ideal $I \subset V^+$
(see B.3) is called {\it friendly} if $(I,V^-)$ is a friendly Jordan pair.
For instance, finite dimensional Jordan pairs over infinite fields are
friendly (because then $U_a$, being the complement of the zero set of the 
polynomial $\det B^\pm(x,a)$, is Zariski-dense in $V^\pm$); 
also, real and complex
Banach Jordan pairs are friendly (since in these cases $U_a$ is an open
neighborhood of the origin), and, more generally, 
all {\it topological (C2)-Jordan pairs} in the sense of [BeNe05, Section 5.1]
are friendly (since also in this case $U_a$ is an 
open neighborhood of the origin). 
The sets $U_a$ essentially describe intersections of
affine parts of $\XX^\pm$, and thus the friendlyness condition assures 
that such intersections are ``not too small". This is
important if one wants to construct some kind of manifold structure
on $\XX^\pm$ (see [BeNe05]).

\Theorem 5.6 (Correspondence between subspaces and sub-Jordan pairs).
Let $(\m,E)$ be a $3$-graded Lie algebra, let 
$V^\pm = \g_{\pm 1}$ be the associated Jordan pair and 
$(\XX^+,\XX^-; \top)$ its associated projective geometry with base
point $(\o^+,\o^-)$. We identify $V^\pm$ with the linear parts 
$(\o^\mp)^\top = \exp(\g_{\pm 1}).\o^\pm$. Then we have the
following correspondence between subspaces and sub-Jordan pairs
(resp., between intrinsic subspaces and inner ideals):
\item{(1)}
If $(\WW^+,\WW^-)$ is a subspace
of the geometry $(\XX^+,\XX^-)$ containing the base point $(o^+,o^-)$,
then  $(W^+,W^-):=(\WW^+ \cap V^+,\WW^- \cap V^-)$ is 
a Jordan sub-pair of the Jordan pair $(V^+,V^-)$, and if
 $\II \subset \XX^+$ is an  intrinsic subspace
containing the base point $o^+$, then
$I := \II \cap V^+$ is an inner ideal  of $V^+$.
\item{(2)}
If $W=(W^+,W^-)$ is a friendly sub-Jordan pair of $(V^+,V^-)$ and
$G_W$ the subgroup of $\PE(\g,E)$
generated by $\exp(W^+)$ and $\exp(W^-)$, 
then $(\WW^+,\WW^-):=(G_W.\o^+,G_W.\o^-)$ is
a subspace of $(\XX^+,\XX^-)$ having $W$ as associated (via $(1)$)
Jordan pair, 
i.e., $W^\pm = \WW^\pm \cap V^\pm$.
If $I$ is a friendly inner ideal of $(V^+,V^-)$, then
$\II := G_I.\o^+$ (where $G_I$ is the group 
generated by $\exp(I)$ and $\exp(V^-)$)
is an intrinsic subspace such that $\II \cap V^+= I$.

\Proof.  (1)
First of all, it follows from the very definition of a subspace
(Section 2.1), that $(W^+,W^-)$ is a pair of {\it linear} subspaces
of $(V^+,V^-)$.
It follows that $W^\pm$ can be naturally identified with its
tangent space $T_{\o^\pm} W^\pm$ at the base point, where the
{\it tangent bundle} $(T\XX^+,T\XX^-)$
can be constructed 
either abstractly (as in [Be02], by using scalar extension
of $\K$ by dual numbers) or directly in the flag model (as in [BeNe04,
Section 2.1]). 
Now, the associated
Jordan pair is constructed in a similar way as the Lie algebra of
a subgroup of a Lie group (see [Be02]): the geometry $(\XX^+,\XX^-)$
(and hence every subgeometry) is a {\it generalized projective geometry}
([Be02, Th.\ 10.1]); for such a geometry, the associated
Jordan pair is nothing but the pair of tangent spaces 
$(V^+,V^-)=(T_{o^+}\XX^+ , T_{o^-}\XX^-)$, equipped with a multilinear
product which is constructed by using the third order tangent bundle
(essentially, this is an algebraic version of the construction of
the associated Lie triple system of a symmetric space, see [Be06,
Chapter 27]). If $(\WW^+,\WW^-)$ is a subgeometry containg
the base point, then the third order tangent bundle 
$(T^3 \WW^+,T^3 \WW^-)$ is a subgeometry of $(T^3\XX^+,T^3 \XX^-)$,
and it follows that the pair of tangent spaces
$(W^+,W^-) \cong (T_{o^+}\WW^+ , T_{o^-}\WW^-)$ is a sub-pair of the Jordan
pair $(V^+,V^-)$.
This proves the first statement,
and the second statement is an immediate consequence since inner ideals
are precisely the subpairs of the form $(I,V^-)$, and intrinsic
subspaces are precisely the subspaces of the form $(\WW^+,\WW^-)$
with $\WW^-$ containing the whole affine part $V^-$.

\ssk (2) Assume  now that a friendly 
sub-pair $(W^+,W^-)$ of $(V^+,V^-)$ is given, and define $\WW^\pm =
G_W.\o^\pm$ as in the claim. The main point of the proof of Part (2) is
to show that
$W^\pm = \WW^\pm \cap V^\pm$.
(Note that this claim means that phenomena such as the ``dense wind"
known from Lie theory do not occur in our situation.)
Now, this claim is a consequence of results from [BeNe04] and [BeNe05]:
in [BeNe04, Th.\ 2.8] it is shown that
the group $\PE(\g,E)$ and hence also $G_W$ act by ``Jordan fractional 
quadratic maps". This means that, whenever a point $p=g.o^+ \in \WW^+$ 
belongs to the affine part $V^+$, 
we can write $g.o^+=d_g(o^+)^{-1} n_g(o^+)$ with a polynomial
$d_g:V^+ \to \End(V^+)$ called ``denominator" and a polynomial
$n_g:V^+ \to V^+$ called ``nominator".
We claim that, if $g \in G_W$, the denominator $d_g(o^+)$ always
belongs to the associative subalgebra $\cal B$ of $\End(V^+)$ generated by
the Bergmann-Operators $B(x,a)$ with $x \in W^+$, $a \in W^-$.
In fact, this is proved exactly as in [BeNe05, Lemma 5.4], by
decomposing $g \in G_W$ into a finite product of elements of the form
$\exp(x_i) \exp(y_i)$ with $x_i \in W^+$, $y_i \in W^-$, 
and by applying induction on the ``length" of such a decomposition.
It is precisely at this point
that the friendlyness condition is needed (cf.\ [BeNe05, loc.\ cit.],
where instead of friendlyness the stronger assumption of a topological
C2-condition was made). 
Having shown that $d_g(\o^+)$ belongs to the algebra $\cal B$ defined
above, it follows that $d_g(\o^+)$ preserves $W^+$ since Bergmann
operators $B(x,a)$ with $x \in W^+, a \in W^-$ preserve $W^+$ because 
$(W^+,W^-)$ is sub-pair.
The same arguments apply as well to the nominators: if $g \in G_W$,
one proves by induction on the decomposition length of $g$ that
$n_g(\o^+)$ belongs to $W^+$.
Finally, we get that
$g.o^+=d_g(o^+)^{-1} n_g(o^+)$ also belongs to $W^+$, i.e.\
$\WW^+ \cap V^+ \subset W^+$. The other inclusion being clear since
$G_W$ contains translations by $W^+$, we have equality. 
Similarly, we see that $\WW^- \cap V^- = W^-$.

\ssk It follows now easily that $(\WW^+,\WW^-)$ is indeed a subgeometry
of $(\XX^+,\XX^-)$. In fact, since, 
by definition, $\WW^+$ and $\WW^-$ are homogeneous under
$G_W$, for every transversal pair $(x,a) \in \WW^+ \times \WW^-$, there
exists $g \in G_W$ such that $g.(x,a)=(\o^+,\o^-)$ (cf.\ [Be02,
Theorem 5.7]). Thus $\WW^+ \cap \v_a$ is a linear subspace of 
$(\v_a,x)$, since the same is true with respect to $\o^-$, as we have just
seen. Similarly, $\WW^- \cap \v_x'$ is a linear subspace of 
$(\v_x',a)$, and the claim on subspaces is proved.
As above, the claim on intrinsic subspaces is an immediate consequence.
\qed

\msk \nin {\bf 5.7. Geometries of intrinsic subspaces.}
In a sense, the preceding theorem shows that the ``geometry of all 
intrinsic subspaces of $\XX^+$" is the global object corresponding
to the ``geometry of inner ideals of $V^+$" (see
Section 5.11 for further comments concerning the latter). 
The geometry of all intrinsic
subspaces being too big, we will again ``fix a type", i.e., we consider
the orbit $G.\II$ of a given intrinsic subspace $\II$  under
the projective group $G=\PE(\g,E)$ and consider the orbit 
as a homogeneous space $G.\II \cong G/G_\II$, where 
$G_\II = \{ g \in G \vert \, g(\II)=\II \}$.
If $\II$ corresponds to the inner ideal $I \subset V^+$, then
(by construction in Theorem 5.6)
$G_\II$ contains the group $G_I = \langle \exp(I),\exp(V^-) \rangle$. 
(It may be strictly bigger: take, e.g., $I=0$; then $G_I = \exp V^-$,
whereas $G_\II $ is the semidirect product of $\exp V^-$ with the
structure group $\Str(V^+,V^-)$, and hence $G.\II$ is the point space
$\XX^+$. In the general case, $G_\II$ is generated by its subgroups
$G_I$ and the linear group $S_I=:\{ g \in \Str(V^+,V^-) \vert \,
g^+(I)=I \}$.) 
As a first step, one may analyse these groups by determining their
``Lie algebras" 
$$
\eqalign{
\g_{I} & = \langle I , \g_{-1} \rangle = I \oplus [I,\g_{-1}] \oplus \g_{-1}, 
\cr
\g_{\II} & = 
I \oplus \s_I \oplus \g_{-1} , \quad \s_I:=
\{ X \in \g_0 \vert \, [X,I] \subset I \}  .\cr}
\eqno (5.1)
$$
The type of intrinsic subspaces we are interested in is defined by the
{\it Peirce-decomposition associated to an idempotent}, cf.\
Appendix A.7 and B.4. 
            
\Theorem 5.8 (Peirce ideals and 5-graded geometries).
Assume $(V^+,V^-)$ is a Jordan pair with corresponding $3$-graded Lie
algebra $\m=\g_1 \oplus \g_0 \oplus \g_{-1}$ and
$(e^+,e^-)$ is a non-zero idempotent in $(V^+,V^-)$, with
associated Peirce-decomposition $V^+ = V_2 \oplus V_1 \oplus V_0$. 
We assume the inner ideal $I:=V_2$ 
to be friendly, and let $\II$ be the  
corresponding intrinsic subspace in $\XX^+$.
\item{(1)} 
There exists a $5$-grading 
$\m=\e_2 \oplus \e_1 \oplus \e_0 \oplus \e_{-1} \oplus \e_{-2}$ 
of $\m$ such that the stabilizer algebra $\g_{\II}$ is equal
to the ``parabolic" subalgebra $\q$ associated to the minus-filtration of this
$5$-grading, i.e.
$$
\g_{\II} = \q := \e_0 \oplus \e_{-1} \oplus \e_{-2}.
$$
\item{(2)}
The group $G_\II$ is equal to the normalizer $N_G(\q)$ of $\q$ in $G$.
In other words, the homogeneous space $G.\II$ is isomorphic to the orbit $G.\q$ 
and also to the $G$-orbit of the minus-filtration of $\e$.

\Proof.
Before constructing the $5$-grading (Lemma 5.9 below), let
us introduce some terminology: 
we say that an inner $5$-grading 
$\m=\h_2 \oplus \h_1 \oplus \h_0 \oplus \h_{-1} \oplus \h_{-2}$ 
of a Lie algebra $\m$ is {\it of Peirce type with respect to
an inner $3$-grading $\m = \g_1 \oplus \g_0 \oplus \g_1$} if
there exist $e^+ \in \g_1$, $e^- \in \g_{-1}$, $H \in \g_0$
such that $\ad(H)$ is the grading derivation of the $5$-grading and
$e^+,H,e^-$ is an $\sL_2$-triple:
$$
[e^+,e^-]=H, \quad [H,e^+]=2e^+, \quad [H,e^-]=-2e^-.
$$

\Lemma 5.9.
Assume $\g=\g_1 \oplus \g_0 \oplus \g_{-1}$ is a $3$-graded Lie algebra
with Euler operator $E$
associated to the Jordan pair $(V^+,V^-)=(\g_1,\g_{-1})$.
\item{(1)} 
Every non-zero idempotent $(e^-,e^+)$ gives rise to a $5$-grading
that is of Peirce type (with respect to the given $3$-grading),
having $H=[e^+,e^-]$ as grading element, and every $5$-grading of Peirce
type is obtained in this way.
\item{(2)}
If $E$ and $H$ are grading elements as above, then $\ad(E)$ and 
$\ad(H)$ are simultaneously diagonalisable, and only the following
combinations of eigenvalues can appear:
$$
\matrix{
\ad(E) \cr \ad(H) \cr} \quad \quad 
\matrix{ -1  & -1 & -1  \cr -2 & -1 & 0 \cr} \quad \quad
\matrix{ 0  &  0 & 0  \cr -1 & 0  & 1 \cr} \quad \quad
\matrix{ 1  &  1 & 1  \cr 0  &  1 & 2 \cr} 
$$
\item{(3)} 
If $E$ and $H$ are grading elements as in $(1)$, then the element
$H':= 2 E - H$ is again a grading element of a $5$-grading $\g =
\e_{-2} \oplus \ldots \oplus \e_2$. 
In terms of the $5$-grading
$\g=\h_2 \oplus \h_1 \oplus \h_0 \oplus \h_{-1} \oplus \h_{-2}$ 
from Part (1) and the $3$-grading 
$\g=\g_1 \oplus \g_0 \oplus \g_{-1}$, we have
$$
\matrix{
\e_2 = \h_0 \cap \g_1 & 
\e_1 = (\h_1 \cap \g_1) \oplus (\g_0 \cap \h_{-1}) &
\e_0 = \h_2 \oplus \h_{-2} \oplus (\g_0 \cap \h_0)    \cr
\e_{-2} = \h_0 \cap \g_{-1} & 
\e_{-1} = (\h_{-1} \cap \g_{-1}) \oplus (\g_0 \cap \h_{1}) &  \cr}
$$
In particular, if $\o=\o(E)$, resp.\ $\f=\f(H')$, 
are the $+$-filtrations corresponding
to the $3$ (resp.\ $5$)-grading just described, we have the relations
$$
\f_2 \subset \o_1 \subset \f_0, \quad \quad \f_1 \subset \o_0 \subset \f_{-1}.
$$
\item{(4)}
Assume that, moreover, the Jordan pair $(V^+,V^-)=(V,V)$ 
is associated to a
Jordan algebra $V$ with unit $1$ and that $(e^+,e^-)=(e,e)$ for an
idempotent $e$ of $V$. Then the grading element $H'$ defined in Part $(3)$
comes again from an idempotent, namely from $(1-e,1-e)$.

\Proof of the lemma. (1)
Assume $(e^+,e^-)$ is a non-zero  idempotent in $(V^+,V^-)$ and let 
$V^+=V_0 \oplus V_1 \oplus V_2$ be the associated Peirce
decomposition. Let $H:=[e^+,e^-] \in \g_0$; then 
we have $[H,e^+]=[[e^+,e^-],e^+]=  T(e^+,e^-,e^+)=
2 Q(e^+)e^- = 2 e^+$ and $[H,e^-]=- 2 e^-$, i.e.,
$e^+,H,e^-$ is an $\sL_2$-triple.
The element $H$ acts diagonally on $V^+$, with possible eigenvalues
$0,1,2$ and the eigenspaces given by the Peirce decomposition
(cf.\ Appendix B, B.4 (1)), and since $H=-[e^-,e^+]$, the same argument
applied to $V^-$ shows that $H$ acts diagonally on $V^-$ with possible
eigenvalues $0,-1,-2$. Since we assume that $\g$ is generated by
$\g_\pm = V^\pm$ and the Euler operator $E$, it follows that
$H$ acts diagonally on
$\g_0 = \sum_{i,j=-2}^2 [(V^+)_i,(V^-)_j] + \K E$, with possible 
eigenvalues $-4,\ldots,4$. But since $V^+$ and $V^-$ are abelian,
the brackets $[(V^+)_i,(V^+)_j]$ etc.\ are zero, and the possible
eigenvalues are among the integers $-2,\ldots,2$. Summing up, the derivation
$\ad(H):\g \to \g$ is diagonalizable with integer eigenvalues 
and defines a 5-grading of $\m$
(cf.\ Prop.\ 5.2).

Conversely, given a $5$-grading of Peirce type, it follows from the
relation $[[e^\pm,e^\mp],e^\pm]=2 e^\pm$ that $(e^+,e^-)$
is an idempotent, whose grading element is, by definition, $[e^+,e^-]$.

\ssk
(2) 
Since $H=T(e^+,e^-)$ belongs to $\g_0$, $\ad(E)$ and $\ad(H)$ commute
and hence are simultaneously diagonalizable. In view of the Peirce
decomposition of $V^+ = \g_1$ and of $V^-=\g_{-1}$ (where only the given
three eigenvalues of $H$ appear), it remains only to prove that
the eigenvalues $2$ and $-2$ of $\ad(H)$ do not appear in $\g_0$.
The proof is merely a question of representation theory of the Lie algebra
$\sL_2$: decomposing $\g$ in isotopycal components under the $\sL_2$-action,
from the given eigenvalues one sees that only the irreducible 
representations of $\sL_2$ on $\rho_j$ on $\K^{j+1}$ for $j=0,1,2$ can
appear, where $\rho_0$ is the trivial representation, $\rho_1$ the natural 
representation by $2 \times 2$-matrices 
$$
\rho_1(H)=\pmatrix{1 & \cr & -1  \cr}, \quad
\rho_1(e^+) = \pmatrix{ 0 & 1 \cr & 0 \cr}, \quad
\rho_1(e^-)= \pmatrix{ 0 & \cr 1 & 0 \cr},
$$
and $\rho_2$ the adjoint representation given by
$$
\rho_2(H)=\pmatrix{2 & & \cr & 0 & \cr & & -2\cr}, \quad
\rho_2(e^+)=\pmatrix{0 & 2 &  \cr &0 & 1 \cr & &0 \cr}, \quad
\rho_2(e^-)=\pmatrix{0 & & \cr 1 &0 & \cr & 2 &0 \cr}.
$$
Now let $\k \subset \g$ be an irreducible $\sL_2$-submodule of $\g$
of type $\rho_2$.
Let $v$ be a highest weight vector: $[e^+,v]=0$, $[H,v]=2v$,
and $v,[e^-,v], [e^-,[e^-,v]]$ spans $\k$ as a $\K$-module.
Thus the lowest weight vector $w= [e^-,[e^-,v]]$ belongs to $V^-=\g_{-1}$,
and the same argument shows that the highest weight vector $v$ belongs to 
$V^+=\g_1$.
Summing up,
the highest weight vector always belongs to $V_2$, and hence 
eigenvectors of $\ad(H)$ for the eigenvalue $2$ always belong
to $\g_1$ and never to $\g_0$. 
In the same way, it is seen that 
eigenvectors of $\ad(H)$ for the eigenvalue $-2$ always belong
to $\g_{-1}$ and never to $\g_0$. 
Since $\g_0$ is $\ad(H)$-stable, we can write $\g_0 = \oplus_{i=-2}^2
(\g_0 \cap \h_i)$; the preceding arguments show that the first and last
term are zero, and hence our claim follows.

\ssk (3) 
Clearly, $\ad(2E-H)$ is again a derivation of $\g$, and using Part (2), 
the eigenvalues are given by the following table: 
$$
\matrix{
\ad(E) \cr \ad(H) \cr \ad(2E-H) \cr} \quad \quad 
\matrix{ -1  & -1 & -1  \cr -2 & -1 & 0 \cr 0 & -1 & -2 \cr} \quad \quad
\matrix{ 0  &  0 & 0  \cr -1 & 0  & 1 \cr 1 & 0  &  -1 \cr} \quad \quad
\matrix{ 1  &  1 & 1  \cr 0  &  1 & 2 \cr 2 &  1 &   0 \cr} 
$$
Thus $\ad(2E-H)$ has again possible eigenvalues $-2,-1,0,1,2$,
and moreover the table gives the eigenspace decomposition stated in the
lemma.
(We call $H'$ the {\it conjugate grading of $H$}. Note that, passing
to the conjugate grading, the eigenspaces belonging to $0$ and $\pm 2$
in $\g_\pm$ are exchanged, whereas those of $\pm 1$ remain unchanged.)
The last relations are also an immediate consequence of eigenvalue
combinations from the table (e.g., we have 
$\f_2 = \e_2 = V_0 \subset V^+ = \g_1 = \o_1$, and so on).

\ssk (4)  
If $e$ is an idempotent in a unital Jordan algebra, then $1-e$ is again
an idempotent, having same Peirce $1$-space, whereas the Peirce $0$-space
and the Peirce $2$-space are exchanged
(cf.\ [Lo75, I.5.6].). As explained above, this is 
exactly what happens in Part (3) for the conjugate grading, and hence
we can conclude that $H'=\ad T(1-e,1-e)$.
\qed

\msk \nin Now we prove Part (1) of the theorem. We define the $5$-grading
$\g=\e_2 \oplus \ldots \oplus \e_{-2}$ as in part (3) of the lemma.
Recall that 
the Lie algebra $\g_{I}$ is generated by $V_2=\g_1 \cap \e_0$
and $V^- = \g_{-1} = \g_{-1} \cap (\e_0 \oplus \e_{-1} \oplus \e_{-2})$,
i.e., 
$$
\g_I = \langle 
(\g_1 \oplus \g_{-1}) \cap (\e_0 \oplus \e_{-1} \oplus \e_{-2}) \rangle.
\eqno (5.2)
$$
Recall further that 
$\s_I = \{ X \in \g_0 \vert \, [X,I] \subset I \}$ 
is the normalizer of $I$ in $\g_0$. We claim that
$$
\s_I = \g_0 \cap (\e_0 \oplus \e_{-1}). 
\eqno (5.3)
$$
The statement (1) follows then from (5.2) and (5.3) since $\g_{\II} = \g_I + \s_I$.
Let us prove (5.3):
the inclusion ``$\supset$" follows directly from the commutation rules
$[\e_i,\e_j] \subset \e_{i+j}$ since $I=V_2 = \e_0 \cap \g_1$. 
In order to prove the converse, assume $X=X_1 + X_0 + X_{-1}$ is the
decomposition of $X \in \g_0$ into eigenvectors of $\ad(H')$.
Assume that $[X,I] \subset I$; this is equivalent to $[X_1,I] \subset I$.
Since, on the other hand, $[X_1,I] \subset \e_1 \cap \g_1 = V_1$, 
our assumption means that $[X_1,I]=0$. We have also $[X_1,V_0] \subset
[X_1,\e_2] \subset \e_3=0$. Therefore, in order to show that $X_1=0$,
all that remains to be proved is that $[X_1,V_1]=0$.
So let $w \in V_1$ and write $w=T(e^+,e^-)w=T(w,e^-)e^+=[T(w,e^-),e^+]$.
Note that $T(w,e^-)=[w,e^-] \in \g_0 \cap \h_{-1}=\g_0 \cap \e_1$.
Then we have
$$
[X_1,w] = [X_1,[T(w,e^+),e^+]]=
[[X_1,T(w,e^+)],e^+] + [T(w,e^+),[X_1,e^+]] = 0
$$
(the first term is zero since $X_1$ and $T(w,e^+)$ belong to the abelian
Lie algebra $\g_0 \cap \h_1$, the second term is zero since $e^+ \in I =
V_2$ and $[X_1,I]=0$ by assumption). It follows that $X_1 =0$, and
thus (5.3) and part (1) are proved.

\ssk Now we deduce part (2) of the theorem.
First of all, note that $\q$ equals its own normalizer, for, if
$X=X_{-2} + \ldots + X_2$ normalizes $\q$, then $[X,H'] \in \q$
since the grading element belongs to $\q$, and it follows that
$2 X_{2} + X_{1} \in \q$, whence $X_{2}=X_{1}=0$.
In a similar way, 
using that, for $X_i \in \e_i$, $e^{\ad(X_i)} H' = H' + [X_i,H'] =
H' - i X_i$,
one sees that the normalizer $N_G(\q)$ is the 
semidirect product of the ``unipotent radical" $\exp(\e_2 \oplus
\e_1)$ and the ``Levi subgroup" $C_G(H')$ (centralizer of the grading
element). The Levi subgroup, in turn, is 3-graded with respect to the
grading element $E$ and hence is generated by $\exp(\e_0 \cap \g_1)$,
$\exp(\e^0 \cap \g_{-1})$ and $C_G(H') \cap C_G(E)$. 
Summing up,
$$
N_G(\q) = \langle \quad \exp(\e_{-2} \oplus \e_{-1} \oplus \e_0 \cap (\g_1 \oplus 
\g_{-1})) , \quad C_G(H') \cap C_G(E) \quad \rangle.
\eqno (5.4)
$$
On the other hand, one deduces from Part (1) that the right hand side of
(5.4) is equal
to $G_\II$. In fact, since $G_\II$ obviously contains the group $G_I$
which acts transitively on $\II$, the main point here is to show that 
for the base point preserving elements we have
$$
G_\II \cap \Str(V^+,V^-) = 
 \langle \quad \exp(\g_0 \cap \e_1) , \quad C_G(H')) \cap C_G(E) \quad \rangle.
\eqno (5.5)
$$
But this follows from (5.3) since $\g_0 \cap \e_1$ is nilpotent and
$C_G(H') \cap C_G(E)$ is precisely the group preserving both gradings.
Summing up, $G_\II = N_G(\q)$.
Finally,  (5.4) shows that $N_G(\q)$ 
is also the group stabilizing the minus-filtration of $\e$, whence the last
statement of Part (2).
\qed

\nin The theorem shows that the geometry of intrinsic subspaces of
Peirce type is isomorphic to some geometry of $5$-filtrations.
Of course, one would then like 
to describe the intrinsic subspace $\II$ directly in terms of filtrations.
By analogy with Theorem 3.8, one is lead to conjecture that $\II$ is the 
set of all $3$-filtrations which satisfy certain inclusion relations
with respect to a fixed $5$-filtration. More precisely:

\msk \nin
{\bf 5.10. Conjecture.} {\it
With notation being as in Theorem 5.8,
the intrinsic subspace $\II$ has the following explicit description in terms
of $5$-filtrations: if $\f=(\f_2 \subset \f_1 \subset \f_0 \subset
\f_{-1} \subset \g)$ is the minus-filtration belonging to the $5$-grading
from Theorem 5.8, then $\II$ is the space of all $3$-filtrations
that are ``squeezed" by $\f$ in the following sense:}
$$
\II = \{ \a \in \XX^+ \vert \, 
\f_2 \subset \a_1 \subset \f_0, \,
\f_1 \subset \a_0 \subset \f_{-1} \}.
$$
 As for the proof of Theorem 3.8, 
 the proof of the inclusion ``$\subset$" is rather elementary, whereas 
the proof of the other inclusion involves a certain transitivity result
which seems  to be related to the problems mentioned at the end of the
introduction and will be taken up elsewhere.


\msk
\nin {\bf 5.11. The linear picture.}
As mentioned in Section 5.7, the space of all inner ideals of $V^+$
can be seen as a sort of ``linear picture" of the geometry of all intrinsic
subspaces of $\XX^+$. 
If we restrict our attention to inner ideals $I$ of 
the type $I=V_2$ for some idempotent $(e^+,e^-)$, this leads us back to
a $3$-graded geometry because 
the Lie algebra $\g_0$ is $3$-graded: with the notation from Lemma 5.9,
$$
\g_0 = (\g_0 \cap \h_1) \oplus  (\g_0 \cap \h_0) \oplus (\g_0 \cap \h_{-1}), 
$$
and hence $(W^+,W^-):=(\g_0 \cap \h_1,\g_0 \cap \h_{-1})$ is a Jordan pair;
let us call it the {\it small Jordan pair associated to an idempotent}.
The linear picture of inner ideals conjugate to $I$ is 
the affine pair geometry corresponding to the small Jordan pair.

\ssk
For the case of certain Banach Jordan triples,
this space is studied
in the paper [Ka01] by W.\ Kaup; more precisely, Kaup considers the space of
{\it complemented principal
inner ideals $I$ in a $JB^*$-triple} (which are all Peirce $2$-spaces for
a suitable idempotent, cf.\ loc.\ cit., Lemma 3.2). 
The principal result in [Ka01] is that this space carries a natural
structure of a complex symmetric Banach manifold (loc.\ cit., Th.\ 4.4),
and that this Banach manifold has as tangent geometry the Peirce 1-space
$V_1$. (The last fact indicates that there is a particular relation between
the small Jordan pair and the Jordan pair $(V_1^+,V_1^-)$;
we do not know if this is a general feature or
if it is due to the particular assumptions in [Ka01].) 
Explicit descriptions are given for the case of Cartan factors:
in case of the series I, II, III, these spaces are again Grassmannians
or products of Grassmannians. 
The additional symmetric space structure from [Ka01] (which does not 
appear in our work) comes from the fact that there Jordan {\it triple
systems} (i.e., $3$-graded Lie algebras with involution, or {\it polar
geometries} in the terminology of [Be02]) are considered,
whereas we work with Jordan pairs, i.e.,  projectively.

\sectionheadline
{Appendix A: Inner ideals in matrix spaces}

\nin {\bf A.1. Left- and right ideals.}
Recall that all left ideals in the associative matrix algebra
$A = M(n,n;\K) = \End(\m)$ (where $\m = \K^n$) 
over a field $\K$ are of the form
$$
I = I_\e = \{ f \in \End(\m) \vert \, \e \subset \ker(f) \}
$$
for some subspace $\e \subset \m$.
(In fact, if $I$ is generated by a single element $g$, then
$I=Ag$, and $f \in Ag$ if and only if $\ker(g) \subset \ker(f)$;
using the relation $I_{\e_1} + I_{\e_2} = I_{\e_1 \cap \e_2}$
and fixing some $\K$-basis $g_1,\ldots,g_r$ 
of $I$, it follows then by induction
that $I = I_{\cap_{i=1}^r \ker(g_i)}$.)
In particular, every left ideal is principal, i.e., generated by a
single element.
By taking transposed matrices, it follows then that all right ideals
in $M(n,n;\K)$ are of the form
$$
J = J_\h = \{ f \in \End(\m) \vert \, \im(f) \subset \h \}
$$
for some subspace $\h \subset \m$.
If $\K$ is merely a ring, then it is in general no longer true that all 
left ideals are of the $\K$-algebra $A$ are of the form $I_\e$.

\msk \nin {\bf A.2. Inner ideals in associative algebras.}
If $A$ is an associative $\K$-algebra, we define a trilinear map
$T:A\times A \times A \to A$ by $T(x,y,z)=xyz+zyx$ and say that
a linear subspace $I \subset A$ is an {\it inner ideal} if it is stable
under multiplication by $A$ ``from the  inside":
$T(I,A,I) \subset I$.
Since $T$ is symmetric in $x$ and $z$, this is equivalent to 
the condition
$Q(I)A \subset I$, where the {\it quadratic map} is defined by
$$
Q(x)y := {1 \over 2} T(x,y,x) = xyx
$$
(by polarization, $T$ may be recovered from $Q$ via
$T(x,\cdot,z)=Q(x+z)-Q(x)-Q(z)$).

\ssk
If $I$ is a left ideal of the associative algebra $A$, then
$Q(I)A \subset IAI \subset II \subset I$, and hence $I$ is
an inner ideal. Similarly, right ideals are inner. Since obviously
the intersection of inner ideals is again an inner ideal, all intersections
of left and right ideals are inner ideals. In particular, if
$A = \End_\K(\m)$ is as above, all subspaces of the form
$$
I_{\e,\h} := I_\e \cap J_\h = \{ f \in \End(\m) \vert \,
\e \subset \ker(f), \, \im(f) \subset \h \}
$$
are inner ideals. 

\msk \nin {\bf A.3. Inner ideals in spaces of rectangular matrices.}
Let
$$
V^+  := M(p,q;\K) = \Hom(\m,\n), \quad \quad
V^-  := M(q,p;\K) = \Hom(\n,\m)
$$
be spaces of rectangular matrices over $\K$, with $\m=\K^q$, $\n=\K^p$.
We define trilinear maps $T^\pm:V^\pm \times V^\mp \times V^\pm
\to V^\pm$ by $T^\pm(x,y,z)=xyz+zyx$ and quadratic maps
$Q^\pm:V^\pm \to \Hom(V^\mp,V^\pm)$ by
$Q^\pm(x)y={1 \over 2} T^\pm(x,y,x)=xyx$.
An {\it inner ideal in $V^+$} is a linear subspace $I \subset V^+$ such
that $Q^+(I)A^- \subset A^+$. As above, it is seen that subspaces
of the form $I_{\e,\h}$ with $\e \subset \m$, $\h \subset \n$, are
inner ideals.

\Proposition A.4.
Assume $\K$ is a field (which, by our permanent assumption, is of
characteristic different from $2$). Then
the inner ideals of $V^+ = \Hom(\m,\n)$ are precisely the sets
$I_{\e,\h} = \{ f \in \Hom(\m,\n) \vert \, \e \subset \ker f, 
\im f \subset \h \}$ for fixed subspaces $\e \subset \m$, $\h \subset \n$.
In other words, with respect to suitable bases in $\m$ and $\n$,
every inner ideal is realized by a space of $p \times q$-matrices of the form
$$
\pmatrix{ A & 0 \cr 0 & 0 \cr}
$$
where $A$ runs through all $r \times s$-matrices, for some fixed
$r \leq q$ and $s \leq p$.
Such an ideal is principal if and and only if $\dim E + \dim F = \dim \m$,
i.e., $r=s$.

\Proof.
This result is well-known in the context of Jordan-theory
(see, e.g., introduction of [Ka01]).
However, we do not know a reference where this result is proved by means
of elementary linear algebra only; 
for the reader who is not specialist in Jordan theory, let us
give such a proof here.

We have already remarked that $I_{\e,\h}$ is indeed an inner ideal.
In order to prove that, conversely, all inner ideals are of this form,
let us first determine all principal inner ideals 
$(g)$, where $g \in \Hom(\m,\n)$.
Since $\K$ is assumed to be a field, we may choose bases such that
$g$ is realized by a matrix of rank $r$ of the form 
$\bigl({\1_r \atop 0}{0 \atop 0}\bigr)$.
An arbitrary element $X \in \Hom(\m,\n)$ is written
accordingly 
$\bigl({a \atop c}{b \atop d}\bigr)$
where $a$ is an $r \times r$-matrix and the other components are 
rectangular. Then we have
$$
Q(g)X=gXg=
\pmatrix{ \1 & 0 \cr 0 & 0 \cr}
\pmatrix{ a & b \cr c & d \cr}
\pmatrix{ \1 & 0 \cr 0 & 0 \cr} = 
\pmatrix{ a & 0 \cr 0 & 0 \cr} .
\eqno {\rm (A.1)}
$$
Thus all such matrices for arbitrary $r \times r$-matrices $a$ do indeed
belong to $(g)$, and it follows that $(g)=I_{\ker(g),\im(g)}$.
This proves that all principal inner ideals are of the form
$I_{\e,\h}$ with $\dim \e + \dim \h = \dim (\m)$.

If $I$ is an arbitrary inner ideal, choose a $\K$-basis
$g_1,\ldots,g_s$ of $I$; then $I$ is the smallest inner ideal
containing the principal ideals $(g_1),\ldots,(g_s)$, and
the claim follows by induction from the following statement:

\Lemma A.5.
Let inner ideals of the form
$I_1:=I_{\e_1,\f_1}$ and $I_2:=I_{\e_2,\f_2}$ be given. Then
the smallest inner ideal $I_1 \lor I_2$ 
containing $I_1$ and $I_2$ is given by
$$
I_{\e_1,\f_1} \lor I_{\e_2,\f_2} =
I_{\e_1 \cap \e_2, \f_1 + \f_2}.
$$
Moreover, we have
$I_{\e_1,\f_1} \lor I_{\e_2,\f_2} = I_{\e_1,\f_1} +
 I_{\e_2,\f_1} +  I_{\e_1,\f_2} +  I_{\e_2,\f_2} $.

\Proof. We fix a basis of $\f_1 \cap \f_2$, complete it to bases of
$\f_1$, resp.\ $\f_2$ and finally to a basis of $\n$.
We order this basis of $\n$ such that first comes the basis
of $\f_1$, then the basis of $\f_2$, with the basis of $\f_1 \cap \f_2$
in the middle, and finally the other basis vectors.
We do the same thing with respect to $\e_1$ and $\e_2$, but we 
order the basis vectors the other way round.
Then $I_1$ is represented by matrices having non-zero
entries on the upper left rectangle, and $I_2$ by matrices having
non-zero entries on a rectangle that possibly intersects the first
rectangle on the basis vectors corresponding to $\f_1 \cap \f_2$ and
$\e_1 \cap \e_2$. 
We may find a quadratic submatrix $g$ in $I_1 + I_2$ having coefficients $1$ on
the diagonal and $0$ else and which is of maximal size with respect
to these conditions. 
The inner ideal condition then implies that all elements
of the form $Q(g)X$ with $X \in V^-$ belong to
$I_1 \lor I_2$, and the matrix calculation (A.1) shows that these
elements fill the whole space of submatrices of the type of $g$.
Iterating this procedure for all possible matrices $g$, we see that
the smallest matrix rectangle containing both $I_1$ and $I_2$ 
belongs to $I_1 \lor I_2$. This proves the lemma and the proposition.
\qed

\Corollary A.6.                      
The rank of a matrix $g \in M(p,q;\K)$ having coefficients in a field
$\K$ is equal to the maximum of lengths of chains 
$0 \subset I_1 \subset \ldots \subset I_r = (g)$
of principal inner ideals contained in the principal ideal $(g)$.

\Proof.
As above, represent $g$ by the matrix
$\bigl({\1_r \atop 0}{0 \atop 0}\bigr)$.
Obviously, 
$$
0 \subset ( \bigl({\1_1 \atop 0}{0 \atop 0}\bigr) )
\subset \ldots \subset
( \bigl({\1_r \atop 0}{0 \atop 0}\bigr) ) =(g)
$$
is a chain of inner ideals, and Prop.\ A.4 shows that it has maximal length.
Its length is $r$, which is the matrix rank of $g$.
\qed

\ssk \nin {\bf A.7. Idempotents.}
Let $(V^+,V^-)$ be a pair of matrix spaces as above.
A pair $(e^+,e^-) \in V^+ \times V^-$ is called an {\it idempotent} if
$$
Q(e^+)e^-=e^+, \, Q(e^-)e^+=e^-, \quad {\rm i.e.} \quad
e^+ e^- e^+ = e^+, \, {\rm and} \, e^- e^+ e^- = e^-.
$$

\Lemma A.8.
A pair $(e^+,e^-) \in \Hom(\m,\n) \times \Hom(\n,\m)$ is an idempotent
if and only if there are bases in $\m$ and $\n$ such that $e^+$ and $e^-$
are simultaneously represented by $p \times q$, resp.\
$q \times p$ matrices of the form
$$
E^+ = \pmatrix{\1_r & 0 \cr 0 & 0 \cr}, \quad
E^- = \pmatrix{\1_r & 0 \cr 0 & 0 \cr}.
$$
In particular, every element $x \in \Hom(\m,\n)$ can be completed to
an idempotent $(e^+,e^-)=(x,y)$.

\Proof.
It is clear that the pair of matrices $(E^+,E^-)$ satisfies the conditions
of an idempotent.
Conversely, assume $(e^+,e^-)$ is an idempotent. 
Decomposing $x \in \m$ in the form
$x=e^-(e^+(x)) + (x-e^-(e^+(x)))$, we see that
$\m = \im(e^-) \oplus \ker(e^+)$, and similarly
$\n=\im(e^+) \oplus \ker(e^-)$. It follows that
$\im(e^-) \to \im(e^+)$, $x \mapsto e^+(x)$ is a bijection with
inverse induced by $e^-$. This map is represented by the identity matrix
$\1_r$ with respect to suitable bases; completing bases by vectors
from the respective kernels leads to the matrix representation $(E^+,E^-)$.

\ssk For the proof of the last statement, represent
$x \in \Hom(\m,\n)$ by the matrix
$E^+$ and let $y \in \Hom(\n,\m)$ be an element represented by the matrix
$E^-$ with respect to the same bases.
\qed

\Proposition A.9 (Peirce decomposition).
If $(e^+,e^-)$ is an idempotent in $V^+ \times V^- =
\Hom(\m,\n) \times \Hom(\n,\m)$,
then the operator
$$
T(e^+,e^-) :  V^+ \to V^+, \quad z \mapsto T(e^+,e^-,z)
$$
is diagonalizable, with three eigenvalues, namely $0$, $1$ and $2$.

\Proof.
We represent $(e^+,e^-)$ by matrices $(E^+,E^-)$ as above. Then
$$
T(E^+,E^-) \pmatrix{a & b \cr c & d \cr} =
E^+ E^-  \pmatrix{a & b \cr c & d \cr} +
 \pmatrix{a & b \cr c & d \cr} E^- E^+ =
   \pmatrix{2a & b \cr c & 0 \cr} .
$$
Hence $T(E^+,E^-)$ is diagonalizable, with eigenspace decomposition
$$
V^+ = V_2 \oplus V_1 \oplus V_0 =
\{ \pmatrix{* & 0 \cr 0 & 0 \cr}  \}
\oplus
\{ \pmatrix{0 & * \cr * & 0 \cr}  \}
\oplus
\{ \pmatrix{0 & 0 \cr 0 & * \cr}  \}.
\qeddis

\nin This description shows that $V_2$ and $V_0$ are inner ideals, according to Section
A.3. Moreover, comparison with Prop.\ A.4 shows:

\Corollary A.10.
Every principal inner ideal of $V^+$ is of the form
$V_2=V_2(e^+,e^-)$ for some idempotent $(e^+,e^-)$.
\qed

\nin
If, in Lemma A.8, we have $r=\min (\dim \m,\dim \n)$,
then $(e^+,e^-)$ is a {\it maximal idempotent}. In this case,
$V_0 = 0$, and then $V_1$ is also an inner ideal.
If  $\dim \m = \dim \n$ and $(e^+,e^-)$ is maximal, then
also $V_1 =0$, i.e.\ $V=V_2$. 
One sees that, if  $\dim \m = \dim \n$,
the collection of all Peirce 2-spaces for suitable idempotents is the same
as the collection of all Peirce 1-spaces.

\msk \nin {\bf A.11. Case of symmetric matrices.}
Assume $\m$ is a $\K$ module together with a strongly non-degenerate
and symmetric bilinear form $\alpha$, i.e., the map
$a:\m \to \m^*$, $x \mapsto \alpha(x,\cdot)$ 
is bijective. Then, for every linear operator
$f:\m \to \m$, the adjoint operator $f^t$ is defined as usual.
The space $V=\Sym(\m)$ of symmetric operators then defines
a pair $(V^+,V^-)=(\Sym(\m),\Sym(\m))$ which is a sub-pair of the Jordan pair
$(\Hom(\m,\m),(\Hom(\m,\m))$. 
Moreover, $V$ is then a {\it Jordan algebra} over
$\K$;  but we will not use this fact here.
{\it Inner ideals} of the pair $(V^+,V^-)$ are defined as above.
The intersection of an inner ideal in the space of all matrices with
the space of symmetric matrices is again an inner ideal; this defines inner
ideals of the form
$$
\{ f \in \Sym(\m) \vert \, \e \subset \ker(f)  \, \im(f) \subset \h \}.
$$
Since  $f$ is symmetric, we have
$\im(f) = \ker(f)^\perp$, and hence the first condition gives us 
$\im(f) \subset \e^\perp$.
Replacing $\e$ by $\e \cap \h^\perp$, we may assume that
$\e = \h^\perp$, and hence the inner ideals in question are all of the  form
$$
I_\e =
\{ f \in \Sym(\o) \vert \, \e \subset \ker(f) \} =
\{ f \in \Sym(\o) \vert \, \im(f) \subset \e^\perp \},
$$
and it is easily seen (same calculation as in (A.1)) that they are all principal.
The proof of the following result is similar to (but more complicated than)
the one of Prop.\ A.4:

\Proposition A.12.
All inner ideals in the Jordan pair $(\Sym(\m),\Sym(\m))$, where 
$\m=\K^n$ is finite-dimensional over a field
$\K$, are principal, and they
are all of the form $I_\e$ for a suitable subspace $\e \subset \K^n$.
\qed

\nin
{\it Idempotents} in $(\Sym(n,\K),\Sym(n,\K))$ are defined as above.
Similarly as in Lemma A.8, one sees that they are just projectors
(i.e., $e^+ = e^- = e$ with $e^2 = e$). In particular, the identity
$\1_n$ is a  maximal idempotent, and as above
(case $\dim \m = \dim \n$),
the collection of all $V_0$'s is the same as the one of all $V_2$'s.

\sectionheadline
{Appendix B: Inner ideals in Jordan pairs}

\msk \nin {\bf B.1. Linear Jordan pairs.}
Let $\K$ be a commutative unital ring in which $2$ and $3$ are invertible.
A pair $(V^+,V^-)$ of $\K$-modules equipped with two trilear maps
$T^\pm:V^\pm \times V^\mp \times V^\pm \to V^\pm$
is called a {\it (linear) Jordan pair} if the following properties are satisfied:

\ssk
\item{(1)}
$T^\pm(x,y,z)=T^\pm(z,x,y)$
\item{(2)}
$T^\pm(x,y,T^\pm(u,v,w))=T^\pm(T^\pm(x,y,u),v,w)-T^\pm(u,T^\mp(y,x,v),w)+              
T(u,v,T(x,y,w))$ 
\ssk

\nin
We use the following notation:
$T^\pm(x,y,x)=2 Q^\pm(x)y$,
$Q^\pm(x,z)=Q^\pm(x+z)- Q^\pm(x)-Q^\pm(z)=T(x,\cdot,z)$.
In any Jordan pair, the {\it fundamental formula}
$$
Q(Q(x)y) = Q(x) Q(y) Q(x)
$$
holds ([Lo75, Prop.\ I.2.2]). For instance, 
associative algebras or spaces of rectangular matrices with $T^\pm$
defined as in Sections A.2 and A.3 are Jordan pairs, and in these
cases the fundamental formula is easily checked by a direct calculation.

\msk \nin {\bf B.2. Construction of Jordan pairs.}
Every linear Jordan pair is obtained as follows: assume 
$\g=\g_1 \oplus \g_0 \oplus \g_{-1}$ is a $3$-graded Lie algebra over
$\K$; we may assume that $\g_0 = [\g_1,\g_{-1}] + \K E$ where $E$
is an Euler operator of the grading (cf.\ Sections 5.1 and 5.5).
Then $(V^+,V^-)=(\g_1,\g_{-1})$ with
$$
T^\pm(x,y,z)=-[[x,y],z]
$$
is a linear Jordan pair, and every linear Jordan pair is obtained in this way
(this is essentially the so-called ``Kantor-Koecher-Tits construction";
cf.\ [Lo95]; see also [Be00, Section III.3] 
for an interpretation of this construction
in terms of polarized Lie triple systems and their standard inbedding).
It is not easy to give a
conceptual interpretation of the fundamental formula in this context; 
on the other hand, it is less difficult to give a conceptual
interpretation of the {\it Bergmann operators}
$$
B^\pm (x,y) := \id_{V^\pm} - T^\pm(x,y,\cdot) + Q^\pm(x) Q^\mp(y)
\, \in \End(V^\pm)
$$
in the context of $3$-graded Lie algebras (see [BeNe04]).
If the pair $(x,y)$ is {\it quasi-invertible} (i.e., $B(x,y)$ is invertible),
then
$$
\beta(x,y) :=(B^+(x,y),B^(y,x)^{-1})
$$
belongs to the {\it structure group} $\Str(V^+,V^-)$ which, by definition,
is the automorphism group of the Jordan pair $(V^+,V^-)$.

\msk \nin {\bf B.3. Inner ideals in linear Jordan pairs.}
A submodule $I \subset V^+$ is called an {\it inner ideal (in $V^+$)} 
if $T^+(I,V^-,I) \subset I$, or, equivalently, $Q^+(I)V^- \subset I$.
We summarize some
general results and notions related to inner ideals
(if not otherwise stated, they are quoted from [Lo75]):

\item{(1)}
For $x \in V^+$, $[x]:= Q(x) V^- = T(x,V^-,x)$ is an inner
ideal, called the {\it principal inner ideal} generated by $x$.
Similarly, for $(x,y) \in V$, $B(x,y) V^+$ is an inner ideal of $V^+$.
Examples: if $(x,y)=(e^+,e^-)$ is an idempotent (see below),
then $V_2 = [e^+]$ is a principal inner ideal and
$V_0 = B(e^+,e^-) V^+$ also is an inner ideal.
\item{(2)}
Let $(x):=[x]+ \K x$, the {\it inner ideal generated by $x$}
(smallest inner ideal containing $x$). Then $x \in [x]$ iff $x$ is 
{\it (von Neumann) regular}, i.e., 
there exists $y \in V^-$ such that $x=Q(x)y$
(in this case $x$ can be completed to an idempotent
$(e^+,e^-)=(x,Q(y)x)$, [Lo75, Lemma I.5.2]).
%
\item{(3)}
An element $z \in V^\pm$ called {\it trivial} if $Q(z)=0$.
Then $\K x$ is an inner ideal, called a  {\it trivial inner ideal}.
$V$ is called {\it non-degenerate} if there are no non-zero trivial elements.
\item{(4)}
A {\it simple} inner ideal $I$ is non-trivial and minimal among non-zero
inner ideals of $V^+$. Equivalently,
$I=V_2^+(d)$ for a division idempotent $d$; or:
$I \not= 0$ and $I=[x]$ for all $0 \not= x \in \m$ (cf.\ [Lo89]).
An element $x \in V$ is called {\it simple} if $(x)$ is a simple inner ideal
(iff $x=d_+$ can be extended to a division idempotent $d=(d^+,d^-)$.
\item{(5)}
{\it descending chain condition} (dcc) on a set $M$ of inner ideals:
every descending chain $\m_1 \supset \m_2 \supset ...$
of inner ideals $\m_i \in M$ becomes stationary.
Similarly: {\it ascending chain condition} (acc).
\item{(6)} images under structural maps:
a {\it structural map} between Jordan pairs $V$ and $W$
is a pair of maps $f:V^+ \to W^+$, $g:W^- \to V^-$ such that
$T^+(f(x),y,f(z))=f T^+(x,g(y),z)$,
$T^-(g(a),b,g(c))=g T^-(a,g(b),c)$.
Then, if $\m \subset V^+$ is an inner ideal, then so is
$f(\m) \subset W^+$.
Example: $W=V^{op}$, $f,g=(Q^+(u),Q^-(u))$;
this gives the principal inner ideals.
\item{(7)} A {\it complement} of an inner ideal $I \subset V^+$ is an inner ideal
$J \subset V^-$ such that $V^+ = I \oplus \ker(J)$ and $V^- = J \oplus \ker(I)$,
where the {\it kernel} is defined by
$\ker(I) = \{ y \in V^- | \, Q^+(I)y = 0 \}$, und dually for $\ker(J)$ (see
[LoNe94]). Note that $\ker(I)$ lives in the space ``dual'' to $I$, and that it is
in general not an inner ideal.
\item{(8)}
annihilators: for $X \subset V^-$  the {\it
annihilator} 
$$
Ann(X):=\Big\{ a \in V^+ | \, \matrix{ Q(a)X=Q(X)a=0  \cr
Q(a)Q(X)=Q(X)Q(a)=0\cr  T(a,X)=T(X,a)=0 \cr} \Big\}.
$$
is an inner ideal. If $e$ is an idempotent,
then $Ann(e^-)=V_0^+(e)$.


\msk \nin
{\bf B.4. Idempotents.}
A pair $e=(e^+,e^-) \in V^+ \times V^-$ such
that $Q(e^+)e^- = e^+$, $Q(e^-)e^+=e^-$ is called
an {\it idempotent}.  
We summarize some basic results and notions related to idempotents:

\item{(1)} Peirce decomposition: the operator $T(e^+,e^-)=T(e^+,e^-,\cdot)$
is diagonalisable with at most three eigenvalues, $0,1,2$; we write
$V = V_0 \oplus V_1 \oplus V_2$ for the eigenspace decomposition of
$V=V^+$, and $V_i=V_i(e)$ is called the {\it Peirce $i$-space} of $e$. 
\item{(2)} Two idempotents 
$e=(e^+,e^-)$, $f=(f^+,f^-)$ called {\it orthogonal} 
if $Q(e^+)f^-=0$ and $Q(e^-)f^+=0$.
Then $e+f=(e^+ + f^+, e^- + f^-)$ is again an idempotent.
\item{(3)}
An {\it orthogonal system} of idempotents is an ordered set of pairwise
orthogonal idempotents; it is called {\it maximal} if it is 
not properly contained in
any larger system. An idempotent $e=(e^+,e^-)$ is called 
{\it maximal} if $\{e \}$ is a maximal orthogonal system.
(Equivalent: the Peirce space $V_0(e)$ contains no non-zero idempotent.)
\item{(4)}
Idempotent $e$ {\it primitive} if it is non-zero and cannot be
written as the sum of two orthogonal idempotents.
\item{(5)}
$e=(e^+,e^-)$ is called a 
{\it division idempotent} $e$ if $V_2(e)=(V_2^+,V_2^-)$ is a division Jordan pair.
(A Jordan pair is called a 
{\it division Jordan pair} if $V \not= 0$ and every non-zero element 
element is {\it invertible}, where $x \in V^\pm$ is called
{\it invertible} if $Q(x): V^\mp \to V^\pm$ is invertible.
A Jordan pair is division iff
it has non-trivial multiplication and no proper inner ideals.)
\item{(6)}
$e=(e^+,e^-)$ is called a
{\it local idempotent} if $V_2(e)$ is a local Jordan par.
(A {\it local Jordan pair} is Jordan pair such that 
the non-invertible elements form
a proper ideal $N$ of $V$; then $V/N$ is division.)
\item{(7)}
{\it Chain condition on idempotents} (cci): no infinite sets of
pairwise orthogonal idempotents (then maximal idempotents exist).
\item{(8)}
A {\it frame} is a maximal orthogonal system of local idempotents.

\def\entries{

\[Be00 Bertram, W., {\it The geometry of Jordan- and Lie structures},
Springer Lecture Notes in Mathematics {\bf vol. 1754},  Berlin 2000

\[Be02 Bertram, W., ``Generalized projective geometries: general
theory and equivalence with Jordan structures", Adv. Geom. {\bf 2}
(2002), 329--369


\[Be04 Bertram, W., ``From Linear Algebra via Affine Algebra to
Projective Algebra", Linear Algebra and its Applications {\bf 378}
 (2004), 109--134

\[Be06 Bertram, W., {\it Differential Geometry, 
Lie Groups and Symmetric Spaces over General Base Fields and Rings.}
Mem.\   AMS, to appear,  arXiv: math.DG/ 0502168 

\[BeNe04 Bertram, W. and K.-H. Neeb, ``Projective completions of
Jordan pairs. I.'' J. Algebra {\bf 277}(2) (2004), 193 -- 225

\[BeNe05 Bertram, W. and K.-H. Neeb, ``Projective completions of
Jordan pairs. II'' Geom.\ Dedicata {\bf 112} (2005), 73 -- 113

\[Bue95 Buekenhout, F. (editor), {\it Handbook of Incidence Geometry},
Elsevier, Amsterdam 1995


\[F73 Faulkner, J.R., ``On the Geometry of Inner Ideals'', J. of Algebra
 {\bf 26} (1973), 1 -- 9

\[F80 Faulkner, J.R., ``Incidence Geometry of Lie Groups in Quantum
Theory'', Group theoretical methods in physics (Proc.\ Eigth Internat.\
Colloq., Kiryat Anavim, 1979) Ann.\ Israel Phys.\ Soc., {\bf 3}, 73 -- 79,
Hilger, Bristol, 1980

\[GiKa98 Gindikin,S.\ and S.\ Kaneyuki, ``On the automorphism group of the
generalized conformal structure of a symmetric $R$-space'',
Diff.\ Geo.\ Appl.\ {\bf 8} (1998), 21 -- 33

\[Hua45 Hua, L.-K., ``Geometries of matrices. I. Generalizations of
von Staudt's theorem'', Trans. A.M.S. {\bf 57} (1945), 441 -- 491

\[Ka01 Kaup, W., ``On Grassmannians associated with $JB^*$-triples",
Math.\ Z.\ {\bf 236}, 567 -- 584, 2001

\[Lo75 Loos, O., {\it Jordan pairs}, Springer LN 460, New York 1975

\[Lo89 Loos, O., ``On the socle of a Jordan pair", Collect.\ Math.\
{\bf 40} (2) (1989), 109 -- 125

\[Lo90 Loos, O., ``Filtered Jordan Systems'', Comm. Alg. {\bf 18} (6)
(1990), 1899--1924

\[Lo91 Loos, O., ``Diagonalization in Jordan Pairs'', J. of Algebra
{\bf 143} (1) (191), 252--268

\[Lo95 Loos, O., ``Elementary Groups and Stability for Jordan Pairs",
K-Theory {\bf 9} (1995), 77 - 116

\[LoNe95 Loos, O.\ and E.\ Neher, ``Complementation of Inner Ideals in 
Jordan Pairs'', J.\ Alg.\ {\bf 166} (1994), 255 -- 295

\[McC04 McCrimmon, K., {\it A Taste of Jordan Algebras}, Springer-Verlag,
New York 2004


}

\references
\dlastpage 

\vfill\eject

\bye